\documentclass[10pt,a4paper]{article}
\usepackage{amsfonts,amsmath,amssymb,accents,amsthm}
\usepackage{verbatim,natbib}
\usepackage{mathrsfs}

\usepackage{xcolor}
\usepackage[colorlinks=true, linkcolor=red, citecolor=black!30!green]{hyperref}

\usepackage[normalem]{ulem}

\usepackage{tikz} 

\newcommand{\dis}{\displaystyle}

\newcommand{\noi}{\noindent}
\newcommand{\halmos}{\rule{1ex}{1.4ex}}
\newcommand{\QED}{\nopagebreak{\hspace*{\fill}$\halmos$\medskip}}

\newcommand{\quand}{\quad\mbox{and}\quad}

\newtheoremstyle{mythm}
  {}
  {}
  {\itshape}
  {}
  {\bfseries}
  {}
  {.5em}
  {#1 #2 \thmnote{(#3)}}

\theoremstyle{mythm}
\newtheorem{theorem}{Theorem}[section]
\newtheorem{proposition}[theorem]{Proposition}
\newtheorem{lemma}[theorem]{Lemma}
\newtheorem{exercise}[theorem]{Exercise}
\newtheorem{corollary}[theorem]{Corollary}
\newtheorem{conjecture}[theorem]{Conjecture}

\newtheorem{counterex}[theorem]{Counterexample}
\newtheorem{remark}[theorem]{Remark}

\newcommand{\bt}{\begin{theorem}}
\newcommand{\et}{\end{theorem}}
\newcommand{\bl}{\begin{lemma}}
\newcommand{\el}{\end{lemma}}
\newcommand{\bp}{\begin{proposition}}
\newcommand{\ep}{\end{proposition}}
\newcommand{\bcor}{\begin{corollary}}
\newcommand{\ecor}{\end{corollary}}
\newcommand{\br}{\begin{remark}\rm}
\newcommand{\er}{\end{remark}}
\newcommand{\bcon}{\begin{conjecture}}
\newcommand{\econ}{\end{conjecture}}
\newcommand{\bex}{\begin{exercise}}
\newcommand{\eex}{\end{exercise}}
\newcommand{\bcou}{\begin{counterex}}
\newcommand{\ecou}{\end{counterex}}

\newenvironment{Proof}[1][]{\noi\textbf{Proof #1}}{\QED}
\newcommand{\bpro}{\begin{Proof}}
\newcommand{\epro}{\end{Proof}}

\newcommand{\be}{\begin{equation}}
\newcommand{\ee}{\end{equation}}
\newcommand{\ba}{\begin{array}}
\newcommand{\ea}{\end{array}}
\newcommand{\bc}{\be\begin{array}{r@{\,}c@{\,}l}}
\newcommand{\ec}{\end{array}\ee}

\newcommand{\al}{\alpha}

\newcommand{\ga}{\gamma}
\newcommand{\Ga}{\Gamma}
\newcommand{\de}{\delta}
\newcommand{\De}{\Delta}
\newcommand{\eps}{\varepsilon}

\newcommand{\ups}{\upsilon}

\newcommand{\Bi}{{\cal B}}

\newcommand{\Gi}{{\cal G}}

\newcommand{\Ri}{{\cal R}}

\newcommand{\R}{{\mathbb R}}
\newcommand{\N}{{\mathbb N}}
\newcommand{\Z}{{\mathbb Z}}

\newcommand{\E}{{\mathbb E}}
\renewcommand{\P}{{\mathbb P}}

\newcommand{\Xb}{{\mathbf X}}

\newcommand{\ti}{\tilde}

\newcommand{\di}{\mathrm{d}}

\newcommand{\pc}{p_{\rm NC}}
\newcommand{\wit}{\widetilde}
\newcommand{\wih}{\widehat}
\newcommand{\xb}{{\bf x}}
\newcommand{\yb}{{\bf y}}
\newcommand{\zb}{{\bf z}}

\newcommand{\sgn}{\text{sgn}}

\newcommand{\tauc}{\tau_{\xb}}

\newcommand{\0}{{\bf 0}}
\newcommand{\mr}{{\rm m}}
\newcommand{\htauc}{\hat\tau_{\xb}}
\newcommand{\hs}{\hspace{1pt}}

\newcommand{\Bl}[1]{B_{#1}\hs l({#1})}
\newcommand{\wb}{\textbf{w}}

\begin{document}


\renewcommand{\labelenumi}{{\rm (\roman{enumi})}}
\renewcommand{\theenumi}{\roman{enumi}}

\author{Jinjiong Yu
	\footnote{KLATASDS-MOE, School of Statistics, East China Normal University, 3663 North Zhongshan Road, Shanghai 200062, China. Email: jjyu@sfs.ecnu.edu.cn}
	\footnote{NYU-ECNU Institute of Mathematical Sciences at NYU Shanghai, Shanghai 200062, China.}
}

\title{Sharp asymptotics for $N$-point correlation functions of coalescing heavy-tailed random walks}

\date{}

\maketitle

\begin{abstract}\noi
We study a system of coalescing continuous-time random walks starting from every site on $\Z$, where the jump increments lie in the domain of attraction of an $\al$-stable distribution with $\al\in(0,1]$.
We establish sharp asymptotics for the $N$-point correlation function of the system.
Our analysis relies on two precise tail estimates for the system density, as well as the non-collision probability of $N$ independent random walks with arbitrary fixed initial configurations. 
In addition, we derive refined estimates for heavy-tailed random walks, which may be of independent interest. 
\end{abstract}
\vspace{.5cm}

\noi


\section{Introduction}


A system of {\em coalescing random walks} (CRW) consists of particles performing i.i.d.\ rate 1 continuous-time random walks, with immediate coalescence when particles meet. 
This fundamental model has attracted sustained interest across disciplines due to its dual role as both a mathematically tractable system and a prototypical example connecting diverse physical phenomena. Its duality with the voter model through graphical constructions \cite{Lig05}, relevance to chemical reaction kinetics \cite{SFH89}, and application to domain wall dynamics in Glauber spin systems \cite{FLM98} underscore its broad significance.

When formulated on $\mathbb{Z}^d$, the CRW can be represented as a $\{0,1\}^{\mathbb{Z}^d}$-valued Markov process, where occupied (resp.\ vacant) sites correspond to 1 (resp.\ 0) states. 
Consider a system of CRW starting from every site on $\Z^d$.
A central quantity of interest is the particle density $\rho_1(t)$, defined via translation invariance as the origin's occupation probability. Pioneering work by Bramson and Griffeath \cite{BG80a}, based on a moment calculation, established dimension-dependent asymptotics for simple random walks:
\begin{equation}
	\rho_1(t) \sim 
	\begin{cases}
		\displaystyle 1/\sqrt{\pi t} & d=1, \\
		\displaystyle \log t/(\pi t) & d=2, \\
		\displaystyle 1/(\varsigma_d t) & d\geq3,
	\end{cases}
\end{equation}
where $\varsigma_d$ denotes the escape probability for $d$-dimensional simple random walks. 
Cox and Perkins \cite{CP04} later provided a probabilistic proof using the super-Brownian invariance principle for the dual voter model. 
For $d\geq3$, where the simple random walk is transient, van den Berg and Kesten~\cite{BK00,BK02} developed a robust approach by deriving an ordinary differential equation for $\rho_1(t)$.

A natural extension of this study is the {\em $N$-point correlation function} $\rho_N(\xb,t)$, defined as the occupation probability at sites $\xb:=(x_1,\ldots,x_N)$ at time $t$, see (\ref{npt2}).
In the marginally recurrent case of coalescing simple random walks in 2D, Lukins, Tribe and Zaboronski \cite{LTZ18} generalized the method of \cite{BK00} and proved
\begin{equation}\label{npt}
	\rho_N(\xb,t)\sim\rho_1(t)^N\times\pc(\xb,t),
\end{equation} 
where the correction term $\pc(\xb,t)$ is the {\em non-collision probability} of $N$ independent walks starting from $\xb$ by time $t$, and it is shown in \cite[Prop~1.3]{CMP10} that $\pc(\xb,t)\sim c_\xb(\log t)^{-n(n-1)/2}$.
This method can be extended to $d\geq3$, but fails for $d=1$ due to the strong correlations induced by the recurrence of the underlying random walk.
Nevertheless, in the case of 1D coalescing Brownian motions, \cite{MRTZ06} derived the asymptotics for $\rho_N(\xb,t)$, relying crucially on the integrability of the system--most notably, the explicit Karlin-McGregor formula.

In contrast to the nearest-neighbor case, much less is known about the $N$-point function for CRW with long-range jumps.
The only existing result in this direction is for 1D coalescing $\alpha$-stable processes with $\al>1$, where the density $\rho_1(t)$ has asymptotic order $t^{1/\al}$, derived as a corollary of \cite[Thm~6.3 and Prop~2.5]{EMS13} and the scaling property of $\al$-stable process.

Motivated by this gap, we study coalescing heavy-tailed random walks on $\Z$, where the jump increments lie in the domain of attraction of $\al$-stable distributions for $0<\al\leq1$. 
Our results reveal that a form of \eqref{npt} remains valid in this setting. Notably, $\al=1$ is the critical point for recurrence/transience, and our study encompasses a class of recurrent random walks. \\ [5pt]
{\noi \bf Notation.} Throughout the paper, we use $c,C$ to denote various positive constants that may even vary within the same line.
For simplicity, when summing over an entire set, we omit the explicit domain and only display the summation variable. For example, for $\xb\in\Z^N$, we write $\sum_{\xb}$ to mean $\sum_{\xb\in\Z^N}$.
For two functions $f(t)$ and $g(t)$, we adopt the convention that
\begin{equation*}\ba{ll}
	\dis f \sim g : \lim_{t\to\infty} f(t)/g(t) = 1, \hspace{1cm}
	&\dis f \asymp g : 0 < \liminf_{t\to\infty} |f(t)/g(t)| \leq \limsup_{t\to\infty} |f(t)/g(t)| < \infty, \\ [5pt]
	\dis f = o(g) : \lim_{t\to\infty} f(t)/g(t) = 0, 
	&\dis f = O(g) : \limsup_{t\to\infty} |f(t)/g(t)| < \infty.
\ea\end{equation*}
To avoid ambiguities, for $k$-th power of a function $f(t)$, we write $f(t)^{k}$ instead of $f^k(t)$. For example, in (\ref{npt}) we do not write $\rho_1^{N}(t)$ but $\rho_1(t)^{N}$.

\subsection{Model and main result}

It is well known \cite{Lig05} that CRW can be constructed by sequentially adding independent random walks, halting each upon its first collision with any of the previously added paths.
Precisely, the coalescing heavy-tailed random walks $(\xi_t)_{t\geq0}$ on $\Z$ can be defined as a $\{0,1\}^\Z$-valued Markov process by specifying its generator thanks to \cite[Chapter I, Thm~3.9]{Lig05}.
For $\xi\in\{0,1\}^\Z$ and $\xb=(x_1,\ldots,x_n)\in\Z^n$, write $\xi(x_1,\ldots,x_n):=(\xi(x_1),\ldots,\xi(x_n))\in\{0,1\}^n$.
We simplify notation by writing $\xi(\xb)=1$ if $\xi(\xb)=(1,\ldots,1)$, and $\xi=1$ if $\xi(x)=1$ for all $x\in\Z$.
Additionally, we represent elements of $\{0,1\}^n$ as binary strings, for example, $(0,1)$ is written as $01$.
With this notation, the generator $\Gi$ of $(\xi_t)_{t\geq0}$ applied on local cylinder functions $f$, i.e., functions that only depend on finitely many coordinates of $\xi$, is given by 
\begin{equation}\ba{l}\label{gen}
	\dis\Gi f(\xi)=\sum_{x,y\in\Z}p(y-x)1_{\{\xi(x,y)=10\}} \big\{f(\xi+e_y-e_x)-f(\xi)\big\}\\ [5pt]
	\hspace{2.2cm}\dis +\sum_{x,y\in\Z}p(y-x)1_{\{\xi(x,y)=11\}}\big\{f(\xi-e_x)-f(\xi)\big\},
	\ea\end{equation}
where $e_x(y)=1_{\{x=y\}}$, and $p(\cdot)$ with $p(0)=0$ and $\sum_x p(x)=1$ denotes the jump kernel, which is in the domain of attraction of an $\al$-stable distribution.
Without loss of generality, we assume that $p(\cdot)$ is aperiodic.
By \cite[Thm~2.6.1]{IL71}, it satisfies
\begin{equation}\label{jump}
	\sum_{|y|\geq x}p(y)=\frac{L(x)}{x^{\al}}\big(1+o(1)\big)\quad\text{as}~x\to\infty,
\end{equation}
where $L:[0,\infty)\to(0,\infty)$ is a {\em slowly varying function} (abbreviated as SVF throughout this paper), meaning that for any $c>0$, it satisfies $L(cx)/L(x)\to 1$ as $x\to\infty$.

To avoid complexity, we assume that $p$ is {\em symmetric} throughout the paper. 
Although, in principle, our results should remain valid in the non-symmetric case, there is a technical obstacle in adapting the present proof, and we therefore leave this extension open (see Subsection~\ref{S:dis} and Remark~\ref{R:nonsym} for further discussion).
For the SVF $L$, there exist $1/\al$-regularly varying functions $b_t:(0,\infty)\to(0,\infty)$ such that $L(b_t)/b_t^{\hs\al}\sim1/t$.
Let $(X_t)_{t\geq0}$ be a rate 1 continuous-time random walk with jump kernel $p$ and initial position $X_0=0$.
It is well known that $X_t/b_t$ converges weakly to an $\al$-stable distribution with probability density function $g_\al(x)$.
By the asymptotic uniqueness of $b_t$ \cite[Thm~1.5.12]{BGT89} and Gnedenko's local limit theorem \cite[Thm~4.2.1]{IL71}, we can choose a specific $1/\al$-regularly varying function
\begin{equation}\label{bt}
	B_t:=\P(X_t=0)^{-1}g_\al(0)
\end{equation}
such that $B_t/t^{1/\al}$ is a SVF,
\begin{equation}\label{btllt}
	\frac{L(B_t)}{B_t^\al}\sim \frac{1}{t}\quand \sup_{x\in\Z}\big|B_t\P(X_t=x)-g_\al(x/B_t)\big|=o(1).
\end{equation}
Recall that $X$ is recurrent if and only if $G(0)=\lim\limits_{x\to\infty}G(0,t)=\infty$, where $G(x,t):=\int_{0}^{t}\P(X_s=x)\di s$ denotes the Green function of $X$ at $x$ up to time $t$.
Following this standard approach, define
\begin{equation}\label{svfl}
	l(t):=\frac{G(0,t)}{g_\al(0)}=\int_{0}^{t}\frac{1}{B_s}\di s, \quad \ell(t):=\frac{B_t\hs l(t)}{t} \quand l(\infty):=\lim\limits_{t\to\infty}l(t)\in(0,\infty].
\end{equation}
In the regime $\al\in(0,1]$,
\begin{equation}\ba{rclcl}
	\dis X~\text{is (marginally) recurrent} & \Leftrightarrow & l(\infty)=\infty, & & \text{if~}\al=1~\text{and}~1/B_t~\text{is not integrable,} \\ [5pt]
	\dis X~\text{is transient} & \Leftrightarrow & l(\infty)<\infty, & & \text{otherwise.}
\ea\end{equation}
In both cases, $l(t)$ is a SVF (hence we adopt the notation $l$ instead of $G$).
Moreover, $\ell(t)$ diverges by Karamata's theorem and is $(\frac{1}{\al}\!-\!1)$-regularly varying.
In the marginally recurrent case, we impose an additional technical assumption to ensure that the recurrence is sufficiently weak, in the sense that the tail probability of the first return time decays slowly.\\ [5pt]
\noi{\bf Function class $\Bi$}: $B_t\in\Bi$ if there exists some $K>0$ such that $\ell(t)^K\geq l(t)$ for all large $t$.\\ [5pt]
Therefore, we focus on two regimes: either the underlying random walk is transient, or it is marginally recurrent with $B_t\in\Bi$.
In our framework, the condition $B_t\in\Bi$ cannot be easily removed due to technical reasons; see Subsection~\ref{S:dis} for further discussion.

To better illustrate $\Bi$, we provide two examples in the marginally recurrent case with $\al=1$.
\begin{itemize}
	\item If $L(t)\sim (\log t)^a$ for some $a<1$, then $B_t\sim t(\log t)^a$, $l(t)\asymp (\log t)^{1-a}$ and $\ell(t)\asymp \log t$, so that $B_t\in\Bi$ with any $K>1-a$.
	\item If $L(t)\sim \exp\{-(\log t)^\kappa\}$ for some $\kappa\in(0,1)$, then $B_t\sim t\hs\exp\{-(\log t)^\kappa\}$, which yields that $l(t)\asymp (\log t)^{1-\kappa}\exp\{(\log t)^\kappa\}$ and $\ell(t)\asymp (\log t)^{1-\kappa}$, and hence $B_t\notin\Bi$.
\end{itemize}

Our main result, Theorem~\ref{T:npf}, concerns the asymptotic behavior of the $N$-point correlation function.
This relies on two key ingredients as in (\ref{npt}).
One is an analysis of the density function $\rho_1(t)$ (Theorem~\ref{T:den}), and the other is an asymptotic result for the non-collision probability of $N$ independent heavy-tailed random walks with fixed starting point $\xb$ (Theorem~\ref{T:noncoa}).
In fact, we go a step further than (\ref{npt}) by providing sharp estimates for these quantities, establishing explicit error bounds in the results below.

\bt[Density asymptotics]\label{T:den} 
Let $(\xi_t)_{t\geq0}$ be the coalescing heavy-tailed random walks on $\Z$ with initial configuration $\xi_0\equiv1$.
Assume the random walk has rate 1 and a symmetric jump kernel $p(\cdot)$ of the form \eqref{jump} with $0<\al\leq 1$.
In the marginally recurrent case, further assume $B_t\in\Bi$ so that $\ell(t)^K\geq l(t)$ for some $K$.
The density of CRW at time $t$ is $$\rho_1(t):=\E[\xi_t(0)]=\P(\xi_t(0)=1).$$ 
Then, there exists some $\eta>0$ such that for any $\eps>0$,
\begin{equation}
	\rho_1(t)=\frac{g_\al(0)\hs l(2t)}{t}
	\big( 1+O\big(\Ri_{1}(t)\big) \big),
\end{equation}
where the error bound $\Ri_1(t)$ can be chosen as follows. 
{\rm (1)} In the marginally recurrent case where $\al=1$, $l(\infty)=\infty$ and $B_t\in\Bi$, we have $\Ri_1(t)=l(t)^{-(\frac{1}{2K}-\eps)}$.
{\rm (2)} In the transient case, if $\al=1$ and $l(\infty)<\infty$, then $\Ri_1(t)=\ell(t)^{-1/2}$; and if $0<\al<1$, then $\Ri_1(t)=t^{-\eta}$.
\et

There exists a natural coupling between coalescing random walks and independent random walks so that the two collections of paths are identical up to the first collision time. For convenience, the following theorem is stated for independent walks.

\bt[Non-collision probability asymptotics]\label{T:noncoa} 
Let $\xb\in\Z^N$, and let $(\Xb_t)_{t\geq0}=(X^1_t,\ldots,X^N_t)_{t\geq0}$, starting from $\Xb_0=\xb$, be $N$ independent continuous-time random walks with rate $1$ and symmetric jump kernel $p(\cdot)$ of the form \eqref{jump} with $0<\al\leq 1$.
In the marginally recurrent case, further assume $B_t\in\Bi$.
Define the first collision time
\begin{equation}\label{xtauc}
	\tauc:=\inf\{t\geq0: X^i_t=X^j_t
	~\text{for some}~1\leq i<j\leq N\}, 
\end{equation}
and set the non-collision probability $\pc(\xb,t):=\P(\tauc>t)$. 
Denote $\ups_N:={N \choose 2}=\frac{N(N-1)}{2}$.
Then, there exists $c_\xb$ depending only on $\xb$ such that
\begin{equation}
	\pc(\xb,t)=\frac{c_\xb}{l(2t)^{\ups_N}}\Big(1+O(\Ri_{\rm NC}(t))\Big),
\end{equation}
where the error bound $\Ri_{\rm NC}(t)=\Ri_{\rm NC}(\xb,t)$ also depends on $\xb$, and can be chosen as follows.\\ 
{\rm (1)} In the case of either $(a)$ marginal recurrence where $\al=1$, $l(\infty)=\infty$ and $B_t\in\Bi$, or $(b)$ marginal transience with $\al=1$ and $l(\infty)<\infty$, we can choose
\begin{equation*}
	\Ri_{\rm NC}(t)=C_\xb\int_{t}^{\infty}\frac{1}{u\hs\ell(u)^{2-\eps}}\di u=o\big(\ell(t)^{-2+\eps}\big)\quad	\text{for any}~\eps>0.
\end{equation*}
{\rm (2)} In the transient case with $0<\al<1$, we can choose $\Ri_{\rm NC}(t)=C_\xb t^{-\eta}$ for any $\eta>1/\al-1$.
\et

We can now present the asymptotic behavior of the $N$-point correlation function of CRW, in line with the spirit of~(\ref{npt}).
\bt[$N$-point correlation function asymptotics]\label{T:npf}
Let $(\xi_t)_{t\geq0}$ be the coalescing heavy-tailed random walks on $\Z$ with initial configuration $\xi_0\equiv1$.
Assume the random walk has jump rate 1 and symmetric jump kernel $p(\cdot)$ of the form \eqref{jump} with $0<\al\leq 1$.
In the marginally recurrent case with $\al=1$ and $l(\infty)=\infty$, further assume that $B_t\in\Bi$.
Denote $\ups_N={N \choose 2}=\frac{N(N-1)}{2}$.
For $N>0$ and $\xb=(x_1,\ldots,x_N)\in\Z^N$, let 
\begin{equation}\label{npt2}
	\rho_N(\xb,t):=\P(\xi_t(\xb)=1)
\end{equation} be the $N$-point correlation function.
Recall $l(t),\ell(t)$ from \eqref{svfl}, and recall $K$ such that $\ell(t)^K\geq l(t)$. 
Recall $\rho_1(t)$ and $\pc(\xb,t)$ from Theorems~\ref{T:den} and \ref{T:noncoa}, respectively.
Then, there exist some $\eta>0$ and constant $c_\xb$ as given in Theorems~\ref{T:noncoa} depending on $\xb$, such that for any $\eps>0$,
\begin{equation}
	\rho_N(\xb,t)=\frac{c_\xb\hs g_\al^N(0)}{t^N\hs l(2t)^{\ups_N-N}} 
	\big( 1+O\big(\Ri_N(t)\big) \big) \sim\rho_1(t)^N\times\pc(\xb,t),
\end{equation}
where the error bound $\Ri_N(t)=\Ri_N(\xb,t)$ also depends on $\xb$, and can be chosen as follows.\\ {\rm (1)} In the marginally recurrent case with $\al=1$, $l(\infty)=\infty$ and $B_t\in\Bi$, for nay $\eps>0$, we can choose $\Ri_N(t)=C_\xb l(t)^{-(\frac{1}{2K}-\eps)}$;\\
{\rm (2)} In the marginally transient case with $\al=1$ and $l(\infty)<\infty$, we have $\Ri_N(t)=C_\xb\ell(t)^{-(\frac{1}{2}-\eps)}$;\\
{\rm (3)} In the transient case with $0<\al<1$, we have $\Ri_N(t)=C_\xb t^{-\eta}$ for any $\eta>1/\al-1$.
\et

\br
The constant $c_\xb$ appearing in Theorems~\ref{T:noncoa} and~\ref{T:npf} depends on both $N$ and the configuration~$\xb$.
It admits an integral representation whose integrand involves the finite-time behavior of the process $(\Xb_t)_{t\ge0}$, see~(\ref{pnc1}) and the discussion above.
However, since this expression is not particularly informative for our purposes, and to keep the statements of the theorems reasonably concise, we omit the explicit formula for $c_{\mathbf{x}}$.
We note, nevertheless, that $c_{\mathbf{x}}$ may vanish as $N \to \infty$.
Therefore, our results do not apply in the regime where one first sends $N \to \infty$ and then $t \to \infty$.
\er

\br
The most interesting case is the recurrent underlying random walk with $\al=1$, $l(\infty)=\infty$ and $B_t\in\Bi$ in Theorems~\ref{T:den} and \ref{T:npf}. 
This mirrors the 2-dimensional simple CRW studied in \cite{BG80a,CP04,LTZ18}.
Indeed, in the 2D setting with our notation one takes $l(t)=\ell(t)=\log t$ and $K=1$, and then $\rho_1(t)$ and $\rho_N(t)$  exhibit the same leading asymptotics as in Theorems~\ref{T:den} and \ref{T:npf}.
Moreover, the error bound in $\rho_1(t)$ can be improved from $O\big(l(t)^{-(1/2-\eps)}\big)$ with any $\eps>0$ to $O\big(l(t)^{-1/2}\big)$.
\er

\subsection{Proof strategy}\label{S:ps}

We first outline the proof of Theorems~\ref{T:den}.
For high-dimensional ($d\geq3$) coalescing random walks, van~der~Berg and Kesten \cite{BK00,BK02} developed a framework to analyze the asymptotics of the system density $\rho_1(t)$ by deriving an ODE through the generator calculations.
This robust method has since been applied for CRW on transitive graphs \cite{HLYZ22} and to coalescing systems with multiple particle types \cite{TZ24}.
Lukins, Tribe, and Zaboronski \cite{LTZ18} adapted this framework to the barely recurrent case in two dimensions.
Furthermore, it is adapted by Lukins, Tribe and Zaboronski \cite{LTZ18} to two-dimensional coalescing random walks which is barely recurrent.
In our setup, we start by applying the generator (\ref{gen}) to the density function to obtain
\begin{equation}\label{rode1}
	\frac{\di}{\di t}\rho_1(t)= -\sum_{y}p(y)\P(\xi_t(0,y)=11)
\end{equation}
To have particles both at sites $0$ and $y$ at time $t$, we can trace back the trajectories for their ancestors (which are typically unique) over a time interval $s$ (with $s\ll t$).
In effect, we consider these as two backward random walks that have not coalesced between $t-s$ and $t$.
By applying the Markov property at time $t-s$, we can decompose
\begin{equation}\label{rode2}
	\dis\P(\xi_t(0,y)=11)\dis\approx \sum_{z_1,z_2}\P(\xi_{t-s}(z_1,z_2)=11) \hspace{1pt} \P(\wih X^{0,y}_{s}=(z_1,z_2),\hat\tau_y>s),
\end{equation}
where $(\wih X^{0,y}_t)_{t\geq0}$ denotes two independent backward random walks starting from $0$ and $y$, respectively, and $\hat\tau_y$ is their first collision time.
Choosing $s=s(t)$ with $1\ll s\ll t$ appropriately ensures that the backward walks are typically far apart at time $s$ (i.e., $|z_1-z_2|\asymp B_{s(t)}$), so that $\P(\xi_{t-s}(z_1,z_2)=11)\approx \P(\xi_{t-s}(z_1)=1)\times\P(\xi_{t-s}(z_2)=1)=\rho_1^2(t-s)$.
Substituting this into (\ref{rode2}) and using $\rho_1(t-s)\approx\rho_1(t)$ (since $t-s$ is close to $t$) yields
\begin{equation}\label{rode3}
	\P(\xi_t(0,y)=11)\approx\rho_1^2(t-s)\P(\hat\tau_y>s)\approx \rho_1^2(t)\P(\hat\tau_y>t),
\end{equation}
where we used the fact that the random walk is either {\em transient} or {\em barely recurrent} (in the sense that $B_t\in\Bi$) so that $\P(\hat\tau_y>s)\approx\P(\hat\tau_y>t)$ for a suitable choice of $s=s(t)$.
Combining (\ref{rode1}) and (\ref{rode3}) produces a simple ODE for $\rho_1(t)$.

Note that (\ref{rode2})-(\ref{rode3}) shows an approximation for the two-point function.
One can extend these ideas to general $\xb\in\Z^N$ to prove Theorem~\ref{T:npf}.
First, we can establish the asymptotic independence
\begin{equation}
	\P(\wih X^{\xb}_s=\zb, \htauc>s)\approx \P(\wih X^{\xb}_s=\zb) \P(\htauc>s),
\end{equation}
since, for either transient or barely recurrent processes, the non-coalescing event $\{\htauc>s\}$ is essentially determined during an initial time interval $[0,u]$ with $u\ll s$, while the positions $\wih X^{\xb}_s$ are mainly determined by the increments from $u$ to $s$. Consequently, we obtain
\begin{equation}
	\rho_N(\xb,t)\approx \sum_\zb\rho_N(\zb,t-s)\P(\wih X^{\xb}_s=\zb, \htauc>s)\approx \rho_1^N(t-s)\hspace{1pt}\P(\htauc>s)\approx \rho_1^N(t)\hspace{1pt}\P(\htauc>t),
\end{equation}
where we use the approximation $\rho_N(\zb,t)\approx \prod_i\rho_1(z_i,t)=\rho_1^N(t)$ for typical well-separated $\zb$.

It remains to establish Theorem~\ref{T:noncoa}, which concerns the non-collision probability $\pc(t)=\pc(\xb,t)=\P(\tauc>t)$ for $N$ independent random walks $(\Xb_t)_{t\geq0}$.
The first step is to derive a crude bound Lemma~\ref{L:pcbdd}, so that we can choose a suitable $s=s(t)\ll t$ satisfying $\pc(s)\approx\pc(t)$.
Then, we derive an ODE for $\pc(t)$ in a manner analogous to the previous arguments, as also done in \cite{LTZ18} for 2D simple CRW. 
Since the decrease in $\pc(t)$ is due to collisions between any two random walks $X^i$ and $X^j$, by the symmetry of $p$ and translation invariance we have
\begin{equation}\label{pc}\begin{array}{r@{\,}c@{\,}l}
	\dis-\frac{\di}{\di t}\pc(t)&=&\dis \sum_{i<j}\sum_{x_{i},x_j}2p(x_i-x_j) \P_{\xb}(\tauc>t,X^i_t=x_{i},X^j_t=x_j) \\ [10pt]
	&=&\dis\sum_{i<j}\sum_{y_{ij}}2p(y_{ij}) \P_{\xb}(\tauc>t,\De^{ij}_t=y_{ij}),
\end{array}
\end{equation}
where $\De^{ij}_t:=X^i_t-X^j_t$ and $y_{ij}:=y_i-y_j$.
Decompose at time $s$ and evolve $\De^{ij}$ backward from $t$ to $s$. 
As only $\De^{ij}_t$ is fixed, other free pairs of walks are unlikely to collide from $s$ to $t$. Thus,
\begin{equation}\label{pc1}
	\P_{\xb}(\tauc>t,\De^{ij}_t=y_{ij}) \approx\sum_{z_{ij}} \P_{\xb}(\De^{ij}_s=z_{ij},\tauc>s) \P_{y_{ij}}(\wih\De^{ij}_{t-s}=z_{ij},\hat\tau_{ij}>t-s).
\end{equation}
Applying asymptotic independence to $\P_{y_{ij}}(\wih\De^{ij}_{t-s}=z_{ij},\hat\tau_{ij}>t-s)$ and noting $|y_{ij}|,|z_{ij}|\ll t^{1/\al}$, the second term in the r.h.s.\ of (\ref{pc1}) can be replaced by $\P_{0}(\wih\De^{ij}_{t-s}=0)\hs \P_{y_{ij}}(\hat\tau_{ij}>t-s)$.
This, together with (\ref{pc}) and a priori estimate $\pc(s)\approx\pc(t)$, implies that
\begin{equation}\label{pc2}
	-\frac{\di}{\di t}\pc(t)\dis \approx \pc(t)\cdot \sum_{i< j}\sum_{y_{ij}}2p(y_{ij}) 
	\P_{y_{ij}}(\hat\tau_{ij}>t-s)\hs \P_{0}(\wih\De^{ij}_{t-s}=0).
\end{equation}
Using standard estimates on the double sum then yields the desired asymptotics for $\pc(t)$.

\subsection{Discussion}\label{S:dis}

To the best of our knowledge, Theorem~\ref{T:npf} is the first result to investigate the $N$-point function for coalescing heavy-tailed random walks, revealing a rich structure of spatial correlations in these systems.
Along the way, we establish several sharp estimates for heavy-tailed random walks (particularly in the regime $\al\leq 1$), which have not appeared in the literature before. We derive quantitative error bounds in asymptotics for key quantities such as the transition probability, the potential kernel, and first hitting times. These estimates are of independent interest and may serve as valuable tools for future investigations into heavy-tailed processes.
Together, these contributions constitute the main novelty of the present work.

In the marginally recurrent case, the assumption that $B_t \in \mathcal{B}$ is of a technical nature.
More precisely, to obtain a quantitative estimate of the non-collision probability $\pc(t)$ in Theorem~\ref{T:noncoa}, an intermediate bound of the type established in Lemma~\ref{L:pcbdd1} is required, where one of the error terms is expressed in terms of $\ell(t)^{-\kappa}$.
If the condition $B_t \in \mathcal{B}$ fails so that $\ell(t)^K \le l(t)$ for all $K$, then in~(\ref{tkl}) the ``error'’ term would, in contrast, dominate the ``main'' term $\pc(t)$ which only involves a power of $l(t)$.
Nevertheless, we expect that the main result should remain valid for $\alpha = 1$ even when $B_t \notin \mathcal{B}$.
But in that case the proof would require an a priori crude bound on $\pc(t)$ in terms of $l(t)$, which would in turn demand new ideas and techniques.
We therefore leave this as an open problem for future work.
We also note a possible extension to the case of non-symmetric jump kernels $p(\cdot)$.
As discussed in Remark~\ref{R:nonsym}, the present proof strategy should still apply provided that Proposition~\ref{P:pk} can be established for the potential kernel in the non-symmetric setting.

Looking ahead, one promising direction is to study the $N$-point function for $\al>1$.
In this regime, the approach of deriving an ODE becomes inapplicable due to the recurrence of the underlying random walk.
One possible avenue worth exploring is to establish a Karlin--McGregor formula for heavy-tailed random walks and to perform an asymptotic analysis of this formula, as is done for coalescing Brownian motions in \cite{MRTZ06}. However, because heavy-tailed paths exhibit crossings caused by jumps, new methods would be required to derive an appropriate Karlin--McGregor formula in this setting.
Another direction is to investigate the scaling limit of the system. By rescaling space by $\eps$ and time by $\eps^\al$ and letting $\eps\downarrow0$, the resulting particle system initially starts from every point of $\R$.
It is known \cite{EMS13} that if $\al>1$, the rescaled system exhibits the phenomenon of coming down from infinity, namely the particles become locally finite at any positive time.
However, this does not occur in our setting where $0<\al\leq 1$. Instead, one may consider the particles as a closed subset of $\R$. We conjecture that, in a suitable topology, the limiting set  (as $\eps\downarrow0$) at a positive time has a fractional Hausdorff dimension of $1-\al$, as suggested by Theorem~\ref{T:npf}. \\ [5pt]
{\noi \bf Organization.} The remainder of the paper is organized as follows. In Section~\ref{S:htrw}, we develop some sharp estimates for heavy-tailed random walks. Section~\ref{S:fct} is devoted to proving Theorem~\ref{T:noncoa} on the non-collision probability $\pc(t)$. Finally, in Section~\ref{S:npt} we establish Theorem~\ref{T:den} for the density function $\rho_1(t)$ and then prove Theorem~\ref{T:npf} concerning the $N$-point correlation function.

\section{Heavy-tailed random walks}\label{S:htrw}

Fundamental results for heavy-tailed random walks, such as transition probabilities, potential kernels, and first hitting times, are essential for our analysis.
Surprisingly, many of these results are not readily available in the existing literature.
Thus, we develop several basic properties in this section.
For future reference, we state our results in the most general setting.

Let $(X^x_t)_{t\geq0}$ be the heavy-tailed continuous-time random walk starting from $X^x_0=x$.
Throughout this section, we assume that the jump rate is 1 and that the symmetric jump kernel $p(\cdot)$ has the form (\ref{jump}).
For brevity, when $x=0$ we also write $(X^0_t)_{t\geq0}$ as $(X_t)_{t\geq0}$.
By \citep[Thm~2.6.5]{IL71} and its proof, the characteristic function of $p$ is given by
\begin{equation}\label{pchf}
	\phi(u)=\exp\big\{-\ga|u|^\alpha\wit L\big(|u|^{-1}\big)\big\},
\end{equation}
where $\wit L$ is also a SVF such that $\wit L\sim L$, and
\begin{equation*}
	\ga=\cos(\frac{\pi\al}{2})\Ga(1-\al)~~\text{for}~\al\neq1, \quand
	\ga=\frac{\pi}{2}~~\text{for}~\al=1,
\end{equation*}
with $\Gamma$ denoting the gamma function.
Let $N_t$ denote the Poisson number of jumps of the random walk by time $t$.
Then the characteristic function of $X_t$ is
\begin{equation}\label{rwchf}
	\phi_t(u)=\E\big[\E[e^{iuX_t}|N_t]\big]=\exp\big(-t(1-\phi(u))\big).
\end{equation}

We will repeatedly use the following two results for a SVF $H(t)$.
The first is the Potter bound \cite[Thm~1.5.6]{BGT89}: There exists $c>0$ such that for any $t\geq s\geq c$ and any $\eps>0$,
\begin{equation}\label{pott}
	(t/s)^{-\eps}\leq H(t)/H(s) \leq (t/s)^{\eps}.
\end{equation} 
The second result is Karamata's theorem \cite[Prop~1.5.8-1.5.10]{BGT89}: For any $c>0$,
\begin{equation}\label{kara}\ba{ll} 
	\text{For}~\al\in(0,1), & \dis\int_t^\infty \frac{H(s)}{s^{1/\al}}\di s \sim \frac{1/\al-1}{t^{1/\al-1}}H(t);\\ [10pt]
	\text{For}~\al=1, & \dis\int_t^\infty \hspace{-2pt}\frac{H(s)}{s}\di s~\text{is a SVF if it is finite, and~} \frac{1}{H(t)}\int_t^\infty\hspace{-2pt}\frac{H(s)}{s}\di s\to0; \\ [10pt]
	& \dis\int_c^t \hspace{-2pt}\frac{H(s)}{s}\di s~\text{is a SVF, and~} \frac{1}{H(t)}\int_c^t\hspace{-2pt}\frac{H(s)}{s}\di s\to0; \\ [10pt] 
	\text{For}~\al\in(1,2), & \dis\int_c^t \hspace{-2pt}\frac{H(s)}{s^{1/\al}}\di s \sim \frac{t^{1-1/\al}}{1-1/\al}H(t).
\ea\end{equation}

\subsection{Transition probability}
\label{S:trprob}
One key input to the ODE framework is a sharp asymptotic estimate of the transition probability $\P(X_t=x)$. 
While Gnedenko’s LLT provides an optimal $o(1)$ error in approximating the transition probability by that of the limiting $\al$-stable process, our analysis requires a finer quantitative error bound. This motivates the following proposition on the difference between transition probabilities at distinct spatial points.
As a supplementary remark, we note that refined local limit theorems with quantitative error bounds are available when the random walk is assumed to lie in the domain of \emph{normal} attraction of an $\alpha$-stable distribution (i.e., the SVF $\widetilde L$ in (\ref{pchf}) is constant and $B_t=t^{1/\al}$), see \cite{CJR23} and references therein.
\begin{proposition}[Transition probability difference]\label{P:dLLT}
	For the heavy-tailed random walk $(X_t)_{t\geq0}$,
	there exists a constant $C$ such that for any $x,y\in\Z$,
	\begin{equation}\label{dLLT}
		\big|\P(X_t=x+y)-\P(X_t=x)\big|\leq \frac{C|y|}{B_{t}^2}.
	\end{equation}
\end{proposition}
\bpro
By the inverse Fourier transform,
\begin{equation*}
	\P(X_t=x)= \frac{1}{2\pi}\int_{-\pi}^{\pi}e^{-iux}\phi_t(u)\di u
	=\frac{1}{2\pi B_t}\int_{-\pi B_t}^{\pi B_t}e^{-ixu/B_t}\phi_t(u/B_t)\di u.
\end{equation*}
Using the inequality $|e^{i(x+y)u/B_t}-e^{ixu/B_t}|\leq |yu|/B_t$,
\begin{equation}\ba{l}\label{trandiff}
	\dis\big|\P(X_t=x+y)-\P(X_t=x)\big|
	\leq \frac{|y|}{2\pi B_t^2} \int_{-\pi B_t}^{\pi B_t} |u\cdot\phi_t(u/B_t)|\di u.
	\ea\end{equation}
Recalling $\phi_t(u)=\exp\{-t(1-\phi(u))\}$ from \eqref{rwchf} with $\phi(u)$ given in (\ref{pchf}), we have
\begin{equation}\label{phit}
	|\phi_t(u/B_t)|= \exp\Big\{-t\Big(1-e^{-\frac{\ga|u|^\al\wit L(B_t/|u|)}{B_t^\al}}\Big)\Big\}
	\leq \exp\Big\{-\frac{c\ga|u|^\al t\wit L(B_t/|u|)}{B_t^\al} \Big\},
\end{equation}
where in the inequality we used that  for small $x>0$, $1-e^{-x}\geq cx$ for some $c>0$. 
By (\ref{btllt}) and the Potter bound (\ref{pott}), for any $0<\eps<\al$,
\begin{equation*}
	\frac{t\wit L(B_t)}{B_t^\al}\frac{\wit L(B_t/|u|)}{\wit L(B_t)} \geq C|u|^{-\eps},
\end{equation*}
which yields
\begin{equation*}
	|u\cdot\phi_t(u/B_t)|\leq |u|\exp\{-t\cdot\frac{C\ga |u|^{\al-\eps}}{t} \} = |u|e^{-C\ga |u|^{\al-\eps}}.
\end{equation*}
Since $|u|e^{-C\ga |u|^{\al-\eps}}$ is integrable over $\R$, the bound (\ref{trandiff}) implies the desired result.
\epro

\br\label{R:nons}
If the random walk has jump rate $r$, one may identify it with a rate 1 random walk evolving over time $rt$. Then,
\begin{equation}\label{dLLTr}
	\big|\P(X_t=x+y)-\P(X_t=x)\big|\leq C|y|/B_{rt}^2 \leq C_r|y|/B_{t}^2.
\end{equation}
The same proof applies when $p$ is non-symmetric. In that case, the characteristic function is
\begin{equation}\label{chfns}
	\phi^{\text{ns}}(u)= \exp\big\{ icu -\ga|u|^\alpha\wit L\big(|u|^{-1}\big)\big(1-i\beta \sgn(u) \omega(u,\alpha) \big) \big\} =\phi(u) e^{i\hs\wit\omega(u)}.
\end{equation}
A slight modification in the bound for $|\phi^{\text{ns}}_t(u/B_t)|$ (using the inequality $|\exp\{-t(1-e^{a+bi})\}|\leq \exp\{-t(1-e^{a})\}$ for $t>0$) yields the same conclusion.
\er

\noi {\bf Large deviation estimates.} For $\al\in(0,2)$, we recall the following large deviation and local large deviation estimates from \cite{Ber19} (which were originally established for discrete-time random walks, but can be adapted to continuous-time via standard Poissonian arguments):
\begin{equation}\ba{cl}\label{ld}
	\text{Large deviation} & \quad \P(|X^0_t| >y)\leq \frac{Ct\1L(y)}{(1+y)^\alpha}\quad{\rm as}~y/B_t\to\infty. \\ [8pt]
	\text{Local large deviation} & \quad\P(X^0_t=y)\leq \frac{Ct\1L(|y|)}{B_t(1+|y|)^\alpha}\quad{\rm for~all}~y\in\Z.
\ea\end{equation}
These bounds are instrumental in controlling the transition probability difference in time.
\bp[Transition probability difference in time]\label{P:dtLLT}
	Let $s\ll t$.
	For the random walk $(X^x_t)_{t\geq0}$, there exists a constant such that
	\begin{equation}\label{dtLLT}
		\dis\big|\P(X^x_t=0)-\P(X^x_{t+s}=0)\big|\leq \frac{C}{B_t}\hs (s/t). 
	\end{equation}
\ep
\bpro
It suffices to show that 
\begin{equation*}
	\big|\frac{\di}{\di t}\P(X^x_t=0)\big|\leq \frac{C}{tB_t},
\end{equation*} 
since the desired result then follows by integrating this inequality from $t$ to $t+s$.
Observe that
\begin{equation*}
	\frac{\di}{\di t}\phi_t(u)=\frac{\di}{\di t} e^{-t(1-\phi(u))} = -(1-\phi(u))\hs e^{-t(1-\phi(u))},
\end{equation*}
Applying the inverse Fourier transform yields
\begin{equation*}
	\Big|\frac{\di}{\di t}\P(X^x_t=0)\Big|
	= \Big|\frac{1}{2\pi}\int_{-\pi}^{\pi} e^{-iux}\frac{\di}{\di t}\phi_t(u) \di u \Big|
	\leq \frac{1}{2\pi}\int_{-\pi}^{\pi} (1-\phi(u))\hs e^{-t(1-\phi(u))}\di u.
\end{equation*}
Using the inequality $x e^{-x} \le e^{-x/2}$ for all $x\ge0$, and setting $x = t(1-\phi(u))$, we obtain
\begin{equation*}
	\Big|\frac{\di}{\di t}\P(X^x_t=0)\Big|\leq t^{-1}\frac{1}{2\pi}\int_{-\pi}^{\pi} e^{-\frac{t}{2}(1-\phi(u))}\di u= t^{-1}\P(X^x_{t/2}=0),
\end{equation*}
where the last equality follows from the identity $e^{-\frac{t}{2}(1-\phi(u))}=\phi_{t/2}(u)$ and the inverse Fourier transform formula.
This concludes the proof by the local limit theorem.
\epro

For non-symmetric kernel $p$, if we replace $X$ by suitably centered random walk, then the large deviation remains valid \cite{Ber19} and thus (\ref{dtLLT}) also holds.

\subsection{Potential kernel}

The potential kernel $a(x)$ of a random walk is defined by
\begin{equation}\label{potker}
	a(x,t):=\int_{0}^{t}\Big(\P(X^0_s=0)-\P(X^x_s=0)\Big)\di s \quand a(x):=\lim\limits_{t\to\infty}a(x,t).
\end{equation}
It is known that for a one-dimensional random walk the limit $a(x)\in[0,\infty)$ exists for all $x\in\Z$ (see \cite[T28.1]{Spi01}.
We need a finer estimate that quantifies the convergence rate.
\begin{proposition}\label{P:pk}
	{\rm(i)} If $0<\al<1$, then for any $0<\kappa_1<\frac{1-\al}{\al}$, we have uniformly in $x$,
	\begin{equation}
		a(x)-a(x,t)= O\big( t^{-\kappa_1} \big).
	\end{equation}
	{\rm(ii)} If $1\leq\al<2$, then there exist $\kappa_1>0$ and $\kappa_2\in(\al/2,1)$ such that
	\begin{equation}
		a(x)-a(x,t)= O\big( t^{-\kappa_1}|x|^{\kappa_2} \big).
	\end{equation}
\end{proposition}
\noi The proof of Proposition~\ref{P:pk} relies on the following elementary integral with small $T$.
\bl\label{L:int}
For any $T>0$ and $x\in\Z$, we have that for all $0\leq \theta\leq 1$,
\begin{equation}
	\int_0^T\frac{1-\cos(xu)}{u^\kappa}\di u=\left\{\ba{ll}
	\dis O\Big(T^{1-\kappa}\Big), &  \text{if~}0<\kappa<1, \\ [10pt]
	\dis O\Big(\log(|x|T)\vee 1\Big), &  \text{if~}\kappa=1, \\ [10pt]
	\dis O\Big(|x|^{\theta+\kappa-1}T^\theta\Big), &  \text{if~}1<\kappa<2.
	\ea\right.
\end{equation}
\el
\bpro
It suffices to treat $x>0$. Using the identity $1-\cos u=\int_{0}^{u}\sin v\di v$, we obtain
\begin{equation}\label{aint}
	\int_0^T\frac{1-\cos(xu)}{u^\kappa}\di u= x^{\kappa-1}\int_0^{xT}\int_{0}^{u}\frac{\sin v}{u^\kappa}\di v\di u= x^{\kappa-1}\int_0^{xT}\int_{v}^{xT}\frac{\sin v}{u^\kappa}\di u\di v.
\end{equation}
It remains to bound the double integral.
Integrating the inner integral, it is easy to see that
\begin{equation*}
	(\ref{aint})\leq \frac{x^{\kappa-1}}{1-\kappa}\int_{0}^{xT}(xT)^{1-\kappa}\sin v \di v \leq C T^{1-\kappa},\quad(0<\kappa<1).
\end{equation*}
\begin{equation*}
	(\ref{aint})\leq \int_{0}^{xT} \sin v\log \frac{xT}{v} \di v \leq C \big(\log (xT)\vee 1\big), \quad(\kappa=1).
\end{equation*}
For $1<\kappa<2$, if $xT\geq1$, then by the integrability of $\sin v/v^{\kappa-1}$ over $(0,\infty)$,
\begin{equation*}
	(\ref{aint})\leq \frac{x^{\kappa-1}}{\kappa-1}\int_{0}^{xT}\frac{\sin v}{v^{\kappa-1}}\di v \leq Cx^{\kappa-1}\leq C x^{\kappa-1}(xT)^\theta, \quad(\theta\geq0).
\end{equation*}
If $xT< 1$, then bounding $\sin v$ by $v$ yields
\begin{equation*}
	(\ref{aint})\leq Cx^{\kappa-1} \int_0^{xT}v^{2-\kappa}\di v\leq Cx^{\kappa-1}(xT)^{3-\kappa} \leq C x^{\kappa-1} (xT)^\theta,\quad(\theta\leq 3-\kappa).
\end{equation*}
This completes the proof.
\epro

\bpro[of Proposition~\ref{P:pk}]
Apply the inverse Fourier transform to $\P(X^0_s=0)-\P(X^x_s=0)$. By Fubini and by the symmetry of $\phi_s(u)=\exp\big(s(1-\phi(u))\big)$, we obtain
\begin{equation}\ba{l}\label{adif}
	\dis a(x)-a(x,t)=\frac{1}{2\pi}\int_{-\pi}^{\pi} \int_{t}^{\infty}\!\! \big(1-e^{ixu}\big) \phi_s(u) \di s\di u =\frac{1}{\pi}\!\int_{0}^{\pi}\!\frac{1-\cos(xu)}{1-\phi(u)} e^{-t\big(1-\phi(u)\big)}\di u.
\ea\end{equation}
Denote $C_0:=\sup \Big\{\ga|u|^\al\wit L\Big(\frac{1}{u}\Big):u\in(0,\pi)\Big\}<\infty$. Using $1-e^{-x}\geq xe^{-x}$, we have
\begin{equation*}
	1-\phi(u)\geq \ga u^\al\wit L\Big(\frac{1}{u}\Big) e^{-C_0}\geq C u^{\al+\eps},\quad (0<u<\pi),
\end{equation*}
where $\eps>0$ can be chosen arbitrary small.
To bound the r.h.s.\ of $(\ref{adif})$, decompose the integral into two parts $D_1:=\{u:0<u\leq t^{-\frac{1-\eps}{\al+\eps}}\}$ and $D_2:=\{u:t^{-\frac{1-\eps}{\al+\eps}}<u\leq\pi\}$. 
On $D_1$,
\begin{equation*}
	I_1:=\int_{D_1} \frac{1-\cos(xu)}{1-\phi(u)} e^{-t\big(1-\phi(u)\big)}\di u
	\leq \frac{1}{C}\int_{0}^{t^{-\frac{1-\eps}{\al+\eps}}} \frac{1-\cos(xu)}{u^{\al+\eps}}\di u.
\end{equation*}
Since $t(1-\phi(u))\geq Ct (t^{-\frac{1-\eps}{\al+\eps}})^{\al+\eps}=Ct^\eps$ on $D_2$,
\begin{equation*}
	I_2\!:=\!\!\int_{D_2}\!\! \frac{1-\cos(xu)}{1-\phi(u)} e^{-t\big(\!1-\phi(u)\!\big)}\di u
	\leq \!\!\int_{D_2}\!\! \frac{1-\cos(xu)}{ C u^{\al+\eps}} e^{-Ct^{\eps}} \di u
	\leq\! \frac{e^{-Ct^{\eps}}}{C}\!\!\int_0^\pi\! \frac{1-\cos(xu)}{u^{\al+\eps}} \di u.
\end{equation*}

For $0<\al<1$, choose $\eps<1-\al$. Then, by Lemma~\ref{L:int}
\begin{equation*}
	a(x)-a(x,t)=I_1+I_2 \leq C t^{-\frac{1-\eps}{\al+\eps}(1-\al-\eps)} +C e^{-Ct^\eps}
	=O\Big(t^{-\kappa_1}\Big),
\end{equation*}
where $\kappa_1:=\frac{1-\eps}{\al-\eps}(1-\al+\eps)$.
Note that $\kappa_1$ approaches $\frac{1-\al}{\al}$ if we send $\eps\to0$.

For $1\leq \al<2$, choose $\eps<\frac{2-\al}{2}$ and apply Lemma~\ref{L:int} with $\theta=2-\al-2\eps$, yielding
\begin{equation*}
	a(x)-a(x,t) \leq C |x|^{1-\eps} t^{-\frac{1-\eps}{\al+\eps}(2-\al-2\eps)} + Ce^{-Ct^\eps}
	=O\Big(t^{-\kappa_1}|x|^{\kappa_2}\Big),
\end{equation*}
where $\kappa_1=\frac{(1-\eps)(2-\al-2\eps)}{\al+\eps}>0$ and $\kappa_2=1-\eps<1$.
\epro

\bl\label{L:ponrel}
For the random walk $(X_t)_{t\geq0}$, its jump kernel $p$ and potential kernel $a$ satisfy
\begin{equation}\label{pkrel}
	\sum_{y}p(y)a(x+y)=1_{\{x=0\}}+a(x).
\end{equation}
\el
\begin{remark}\label{R:jrr}
	This is analogous to the relation between the Green's function and the jump kernel of a transient random walk.  
	We will use \eqref{pkrel} for $x$=0.
	Observe that if the random walk is of rate $r$, then $\sum_{y}p(y)a(y)=r$.
\end{remark}
\bpro
Consider the event $\{X_t^0=0\}$ and decompose it according to whether the walk makes at least one jump during the time interval. Let $N_{s,t}$ denote the number of jumps during $(s,t]$.
\begin{equation*}\ba{l}
	\dis \P(X_t^x=0)=\int_{0}^{t}\sum_{y\neq 0}\P(X^x_s=-y)p(y)\P(N_{s,t}=0) \di s + \P(N_{0,t}=0)\cdot 1_{\{x=0\}} \\ [5pt]
	\dis \hspace{1.88cm} = \int_{0}^{t}\sum_{y\neq 0}\P(X^{x+y}_s=0)p(y)e^{s-t}\di s + e^{-t}1_{\{x=0\}}.
	\ea\end{equation*}
Integrating in $t$ from 0 to $T$ in the above equation and using Fubini, we have
\begin{equation*}\ba{l}
	\dis\int_{0}^{T}\P(X_t^x=0)\di t
	= \int_{0}^{T} \sum_{y\neq 0}\P(X^{x+y}_s=0)p(y) (1-e^{s-T})\di s + (1-e^{-T})1_{\{x=0\}}.
	\ea\end{equation*}
Moving the term $\int_{0}^{t} \sum_{y\neq0} \P(X^{x+y}_s=0)p(x)\di s$ to the r.h.s.\ and by the definition of $a(x,T)$,
\begin{equation*}\ba{l}
	\dis\sum_{y}p(y)a(x+y,T)-a(x)= 
	1_{\{x=0\}} -e^{-T} \Big(\int_{0}^{T}\sum_{y\neq 0}\P(X^{x+y}_s=0)p(y)e^{s}\di s +1_{\{x=0\}}\Big).
\ea\end{equation*}
Since $\sum_{y\neq 0}\P(X^{x+y}_s=0)p(x)\leq \max_x\P(X^{x+y}_s=0)\to0$ as $s\to0$, letting $T\to\infty$ yields the desired result.
\epro

\br\label{R:nonsym}
In the non-symmetric case, where the characteristic function is replaced by $\phi^{\mathrm{ns}}$ in~(\ref{chfns}), an equation analogous to~(\ref{adif}) holds for the symmetrized quantity $a(x) - a(x,t) + a(-x) - a(-x,t)$.
Hence, Proposition~\ref{P:pk} extends to this symmetrized form.
However, deriving the result for $a(x) - a(x,t)$ alone presents certain technical difficulties that may require a more refined analysis, particularly in the case $\alpha = 1$ when the mean does not exist; see, for instance, \cite[Section~4]{Ber19} for further discussion.
We note that, in our arguments, the required estimates for a single random walk rely solely on the results from this subsection and the previous one (Remark~\ref{R:nons}). The remaining estimates, particularly those involving collision times, can be reduced to the first hitting times of difference process between two random walks. Notably, the jump kernel of the difference process is always symmetric.
Consequently, the main result, Theorem~\ref{T:npf}, would continue to hold provided Proposition~\ref{P:pk} can be extended to non-symmetric jump kernels~$p(\cdot)$.
\er

\subsection{First hitting time}\label{S:hit}

For the heavy-tailed random walk $(X^x_t)_{t\geq0}$, the first hitting time of 0 is defined by
\begin{equation}\label{taux}
	\tau_{0,x}:=\tau_x:=\inf\{t:X_t^x=0\text{~and~}X^x_s\neq0~\text{ for some~} 0<s<t\}.
\end{equation}
Our goal is to derive asymptotic estimates for the tail probability $	\P(\tau_x>t)$.
A classic approach (which we leave the details to readers) is to apply the Laplace transform to the last-exit formula and then invoke a Tauberian theorem \cite[Theorem~1.7.1]{BGT89}, which yields
\begin{equation}\label{hit0}
	\P(\tau_x>t) \sim \frac{a(x)}{c_\al g_\al(0)l(t)},
\end{equation}
where $c_\al=1$ for $0<\al\leq1$ and $c_\al=\Ga(1/\al)\Ga(2-1/\al)$ for $1<\al<2$.
We leave the details to readers.
Nevertheless, to incorporate these estimates into our ODE framework, we require a higher order estimate of the tail probability. 
Fortunately, for $\al\leq1$ the function $\P(\tau_x>t)$ is slowly varying, which allows us to obtain the refined result Proposition~\ref{L:taut}.

Below, we only focus on case $0<\al\leq1$. 
Define the auxiliary function
\begin{equation}\label{jt}
	j(t):=\sum_{x\in\Z}p(x)\P(\tau_x>t)=\sum_{x\neq0}p(x)\P(\tau_x>t),
\end{equation}
which represents that the random walk jumps away from 0 initially and does not hit 0 by time $t$.
Recall the last-exit formulae that decompose at the last time the walk is at $0$,
\begin{equation}\ba{l}\label{exit}
	\dis 1=\P(X^0_t=0)+\int_{0}^{t}\P(X^0_s=0)j(t-s)\di s, \\ [5pt]
	\dis \P(\tau_x\leq t)=\P(X^x_t=0)+\int_{0}^{t}\P(X^x_s=0)j(t-s)\di s, \quad(x\neq0),
	\ea\end{equation}
We begin by estimating $j(t)$.
Recall that $g_\al$ is the density function of the limiting $\al$-stable distribution, and recall from (\ref{svfl}) that $l(t)=\int_0^t \frac{1}{B_s}\di s$ and $\ell(t)=B_t\hs l(t)/t$.
\bl\label{L:jest}
For $j(t)$ defined in \eqref{jt} and for $\al\leq1$, one has
\begin{equation}\label{jest}
	\sum_xp(x)\P(\tau_x>t) = j(t)=\frac{1}{g_\al(0)l(t)}\Big(1+O(1/\ell(t))\Big).
\end{equation}
\el
\bpro
Since $j(t)=\sum_xp(x)\P(\tau_x>t)$ is decreasing in $t$, the last-exit formula (\ref{exit}) implies
\begin{equation*}
	1\geq\P(X_t^0=0)+j(t)\int_0^t \P(X_s^0=0)\di s\geq g_\al(0)l(t)j(t).
\end{equation*}
For the other inequality, we replace $t$ by $2t$, use $\P(X^0_s=0)\leq C/B_t$ for $s\in[t,2t]$, and obtain
\begin{equation*}\ba{l}
	\dis 1\leq \frac{C}{B_{t}} + \int_t^{2t}\frac{C}{B_t}j(2t-s)\di s  +j(t)\int_{0}^{t}\P(X_s^0=0)\di s \leq \frac{C}{B_{t}} + \frac{Ct}{B_tl(t)}+g_0(\al)l(t)j(t),
\ea\end{equation*}
where in the last inequality, for the integral from $t$ to $2t$ we used $j(t)\leq \frac{1}{g_\al(0)l(t)}$ and then $\int_0^t\frac{1}{l(s)}\di s\sim\frac{t}{l(t)}$.
Note that $B_t\gg\ell(t)$. The estimate (\ref{jest}) follows.
\epro

\noi We now turn to the difference of $j$ over time.
\bl\label{L:dj}
For $j(t)$ defined in \eqref{jt} and for $\al\leq1$,  there exists $C$ such that for $s<t$ and any $\eps>0$,
\begin{equation}\label{dj}
	j(s)-j(t)\leq \frac{C}{l(t)\ell(t)} \left(\frac{t}{s}\right)^{\frac{1}{\al} - 1 + \eps}.
\end{equation}
\el
\bpro
Since $t^{1/\al}/B_t$ is a SVF, the Potter bound implies that for any $\eps>0$,
\begin{equation}\label{dl}\ba{l}
	\dis l(t)-l(s)=\frac{t^{1/\al}}{B_t}\int_{s}^{t}\frac{u^{1/\al}/B_u}{t^{1/\al}/B_t} \frac{1}{u^{1/\al}}\di u \leq \frac{t^{1/\al}}{B_t}\int_{s}^{t}\frac{(t/u)^\eps}{u^{1/\al}}\di u  \\ [10pt]
	\hspace{5cm}\dis\leq \frac{t}{B_t}\left(\frac{t}{s}\right)^{\frac{1}{\al} - 1 + \eps} 
	= \frac{l(t)}{\ell(t)}\left(\frac{t}{s}\right)^{\frac{1}{\al} - 1 + \eps}.
\ea\end{equation}
Similarly, for the SVF $l(t)$ and $(1/\al-1)$-regularly varying function $\ell(t)$,
\begin{equation*}
	1/l(s)\leq (t/s)^\eps/l(t),\qquad
	1/\ell(s) \leq (t/s)^{1/\al-1+\eps}/\ell(t).
\end{equation*} 
Inserting these bounds into (\ref{jest}) yields that for any $\eps>0$,
\begin{equation*}
	j(s)-j(t)\leq \frac{l(t)-l(s)}{g_\al(0)l(t)\hs l(s)} +\frac{C}{l(t)\ell(t)} +\frac{C}{l(s)\ell(s)}  \leq \frac{C}{l(t)\ell(t)}  (t/s)^{\frac{1}{\al}-1+\eps},
\end{equation*}
which completes the proof.
\epro

We are now ready to deduce the tail probability $\P(\tau_x>t)$.
For convenience, denote 
\begin{equation}
	\De(x,t):=\P(X_t^0=0)-\P(X_t^x=0)\geq0,
\end{equation} 
so that the last-exit formula \eqref{exit} yields, for $x\neq0$,
\begin{equation}\label{xexit}
	\P(\tau_x>t)=\De(x,t)+\int_0^t\De(x,s)j(t-s)\di s.
\end{equation}
Note that a uniform bound on the transition probability $\P(X^x_t=0)\leq C/B_t$ implies the uniform bound $\Delta(x,t)\leq C/B_t$.
This, combined with Proposition~\ref{P:dLLT}, shows that
\begin{equation}\label{Debdd}
	0\leq \De(x,t) \leq CB_t^{-1} (|x|B_t^{-1}\wedge 1). 
\end{equation}

\bp[Higher order estimate for hitting probability]\label{L:taut}
Let $(X^x_t)_{t\geq0}$ be the heavy-tailed random walk with increments in the domain of attraction of an $\al$-stable law for some $0<\al\leq1$.
Recall the first hitting time $\tau_x$ from \eqref{taux}. Recall $B_t$ from \eqref{bt}, and recall $l(t)$ and $\ell(t)$ from \eqref{svfl}.
Then, for $x\neq0$ and any $\eps>0$, there exists some $\kappa>0$ such that
\begin{equation}\label{taut}
	\P(\tau_x>t) = \left\{\ba{ll}
	\dis \frac{a(x)}{g_\al(0)l(t)} + O\Big(\frac{a(x)}{l(t)\ell(t)} + \frac{|x|}{t^{\kappa}}\Big), & \al=1; \\ [12pt]
	\dis \frac{a(x)}{g_\al(0)l(t)} + O\Big(\frac{a(x)}{l(t)\ell(t)} + \frac{t^\eps}{\ell(t)}\Big), & 0<\al<1.
\ea\right.\end{equation}
\ep
\bpro
Rewrite (\ref{xexit}) as
\begin{equation}\label{taux1}
	\P(\tau_x>t)= \int_{0}^{t}\De(x,s)j(t)\di s+ \int_{0}^{t}\De(x,s)\big(j(t-s)-j(t)\big)\di s + \De(x,t).
\end{equation}
We then estimate the r.h.s. Firstly, by Proposition~\ref{P:pk} and Lemma~\ref{L:jest},
\begin{equation}\label{taubdd1}
	\int_{0}^{t}\De(x,s) j(t)\di s= a(x)j(t)+\big(a(x,t)-a(x)\big)j(t)= \frac{a(x)}{g_\al(0)l(t)} + O\Big( \frac{a(x)}{l(t)\ell(t)}+ \frac{|x|^{\kappa_2}}{t^{\kappa_1}l(t)}\Big).
\end{equation}
In the equation above, if $\al=1$, then $\kappa_1>0,\kappa_2\in(0,1)$ and thus $\frac{|x|^{\kappa_2}}{t^{\kappa_1}l(t)}\leq C\frac{|x|}{t^\kappa}$ for $\kappa\in(0,\kappa_1)$; and if $0<\al<1$, then $\kappa_1\in(0,\frac{1-\al}{\al}),\kappa_2=0$ and thus 
 
where $\kappa_1>0$ and $\kappa_2\in(0,1)$ if $\al=1$, and $\kappa_1\in(0,\frac{1-\al}{\al})$ and $\kappa_2=0$ if $0<\al<1$.
Secondly, applying Lemma~\ref{L:dj} to $j(t-s)-j(t)$ for $s\in(0,t/2)$,
\begin{equation}\label{taubdd2}
	\int_{0}^{\frac{t}{2}}\De(x,s)\big(j(t-s)-j(t)\big)\di s \leq \frac{C}{l(t)\ell(t)}\int_{0}^{\infty}\De(x,s) \di s = \frac{Ca(x)}{l(t)\ell(t)}.
\end{equation}
For $s\in(t/2/,t)$, Karamata's theorem implies $\int_0^{t}j(s)\di s \asymp t/l(t)$, and thus by (\ref{Debdd})
\begin{equation}\label{taubdd3}
	\int^{t}_{\frac{t}{2}}\De(x,s)\big(j(t-s)-j(t)\big)\di s \leq \frac{C}{B_{t}}\Big(\frac{|x|}{B_t}\wedge1\Big)\int_{0}^{\frac{t}{2}} \big( j(s) +j(t) \big) \di s 
	\leq \frac{C}{\ell(t)}\Big(\frac{|x|}{B_t}\wedge1\Big).
\end{equation}
Recall that $B_t\geq\ell(t)$ for large $t$. (\ref{Debdd}) implies that $\De(x,t)\leq C\ell(t)^{-1}(|x|B_t^{-1}\wedge 1)$.
 
For $\al=1$, since both $l$ and $\ell$ are SVFs, there exists $\kappa\in(0,1)$ such that both $\frac{|x|^{\kappa_2}}{t^{\kappa_1}l(t)}$ and $\frac{|x|}{\ell(t)B_t}$ are bounded by $|x|/t^{\kappa}$.
Using this bound for the r.h.s\ of (\ref{taubdd1}) and (\ref{taubdd3}), and inserting (\ref{taubdd1})-(\ref{taubdd3}) into (\ref{taux1}), we obtain (\ref{taut}) for $\al=1$.
For $0<\al<1$ and for any $\eps>0$ small enough, recall $\kappa_2=0$ and choose $\kappa_1=1/\al-1-2\eps$ in (\ref{taubdd1}). Since now $l$ is a SVF and $\ell$ is $(1/\al\!-\!1)$-regularly varying, $\frac{1}{t^{\kappa_1 }l(t)}$ is bounded by $t^\eps/\ell(t)$. 
Thus, inserting (\ref{taubdd1})-(\ref{taubdd3}) into (\ref{taux1}), we obtain (\ref{taut}) for $0<\al<1$.
\epro

Since $\ell(t)=B_t\hs l(t)/t$, we observe that $B_t\gg\ell(t)^{\theta}$ if $0<\theta<\frac{1}{1-\al}$ (with the convention that $\frac{1}{1-\al}=\infty$ if $\al=1$).
The following lemma gives a tail bound on the hitting time when the starting position is relatively far from the origin.
\bl\label{L:tau<}
Let $\theta\in(0,\frac{1}{1-\al})$.
Under the assumptions of Proposition~\ref{L:taut}, if $|x|\geq B_t/\ell(t)^\theta$, then for any $\eps>0$ there exists $C>0$ such that
\begin{equation}\label{tau<}
	\P(\tau_x\leq t)\leq C\hs \ell(t)^{-1+(1-\al)\theta+\eps}.
\end{equation}
\el
\bpro
We use the last-exit formula (\ref{exit}).
Recall from (\ref{btllt}) that $L(B_t)/B_t^\al\sim1/t$. By the Potter bound
\begin{equation*}
	\frac{L\big(B_t/\ell(t)^{\theta}\big)}{B_t^\al/\ell(t)^{\al\theta}}
	\leq \frac{C \ell(t)^{\al\theta}}{t}\cdot \frac{L\big(B_t/\ell(t)^{\theta}\big)}{L(B_t)} \leq \frac{C\ell(t)^{\al\theta+\eps} }{t}.
\end{equation*}
Set $s:=t/\ell(t)^{\al\theta}$.
Since the polynomial power in $\ell(t)^{\al\theta}$ is $\frac{1-\al}{\al}\cdot\al\theta<1$, we have $s\gg1$.
We now split the integration in (\ref{exit}) into three regions.
For $u\in(0,s]$, by the local large deviation (\ref{ld}), Karamata's theorem, and $j(t-u)\leq C/l(t)$,
\begin{equation*}
	\int_0^{s} \P(X^x_u=0)j(t-u)\di u\leq \int_0^{ s}\frac{CuL\big(B_t/\ell(t)^{\theta}\big)}{B_uB_t^\al/\ell(t)^{\al\theta}} \cdot\frac{1}{l(t)}\di u \leq \frac{C s^2\ell(t)^{\al\theta+\eps}}{t B_{ s}\hs l(t)} \leq \frac{C \hs (t/ s)^{\frac{1}{\al}-2+\eps}}{\ell(t)^{1-\al\theta-\eps}}.
\end{equation*} 
For $u\in( s,t- s]$, we bound $j(t-u)\leq j( s)\leq C/l( s)$, and use (\ref{dl}) to obtain 
\begin{equation*}
	\int^{t- s}_{ s} \P(X^x_u=0)j(t-u)\di u \leq \frac{C}{l( s)}\big( l(t)- l( s)\big)\leq \frac{C}{\ell(t)}(t/ s)^{\frac{1}{\al}-1+\eps}.
\end{equation*}
Finally, for $u\in(t- s,t]$, by the uniform bound $ P(X^x_u=0)\leq C/B_t$ and Karamata's theorem,
\begin{equation*}
	\int_{t- s}^{t} \P(X^x_u=0)j(t-u)\di u \leq \frac{C}{B_t}\int_0^{ s}j(u)\di u \leq \frac{C s}{B_tl( s)}\leq \frac{C}{\ell(t)}(t/ s)^{-1+\eps}.
\end{equation*}
Recalling that $B_t$ is $1/\al$-varying and $\ell(t)$ is $(1/\al-1)$-varying, we have $1/B_t\leq C/\ell(t)$. Thus,
\begin{equation*}
	\P(X^x_t=0)\leq C/B_t\leq C/\ell(t).
\end{equation*}
Combining the four bounds above and inserting $t/s=\ell(t)^{\al\theta}$ gives
\begin{equation*}
	\P(\tau_x\leq t)=\P(X^x_t=0)+\int_{0}^{t}\P(X^x_u=0)j(t-u)\di u 
	\leq \frac{C}{\ell(t)^{1-(1-\al)\theta-\eps}}.
\end{equation*}
where $\eps>0$ can be chosen arbitrary small.
\epro

We also include a lemma dealing with the joint event $\P(\tau_x>t,X^x_t=y)$.
\bl\label{L:tauXt}
Under the assumptions of Proposition~\ref{L:taut},
\begin{equation}\label{tauxt}
	\P(\tau_x>t,X^x_t=y)\leq \frac{Ca(x)}{B_t\hs l(t)} \quand \sum_xp(x)\P(\tau_x>t,X^x_t=y)\leq \frac{C}{B_t\hs l(t)}.
\end{equation}
Moreover, for any $\kappa< \al/(1-\al)$, $($if $\al=1$ then $\kappa<\infty$$)$,
\begin{equation}\label{ptauXt}
	\sum_xp(x)\Big|\P(\tau_x>t,X^x_t=y)-\P(\tau_x>t)\P(X^0_t=0)\Big| \leq \frac{C }{B_tl(t)\ell(t)^{\kappa}}+\frac{C|y|}{B_t^2}.
\end{equation}
\el
\bpro
For the bound in (\ref{tauxt}), we decompose at time $t/2$, by the Markov property,
\begin{equation*}
	\P(\tau_x>t,X^x_t=y)\leq 
	\sum_z\P(\tau_x>t,X^x_t=y)\P(X^z_{t/2}=y) \leq\frac{C}{B_{t/2}}\P(\tau_x>t/2) \leq\frac{Ca(x)}{B_tl(t)}.
\end{equation*}
The second bound follows upon summing over $x$ and applying Lemma~\ref{L:ponrel}.

For the estimate in (\ref{ptauXt}), note that if $|x|>B_t^2$, then for any  $\eps>0$,
\begin{equation}
	\sum_{|z|>B_t^2}p(x)\P(\tau_x>t,X^x_{t}=y)\leq \sum_{|z|>B_t^2}p(x)\leq \frac{CL(B_t^2)}{B_t^2}\leq \frac{C}{B_t^{2-\eps}}.
\end{equation}
It remains to treat $|x|\leq B_t^2$.
Let $H_1=H_1(t), H_2=H_2(t)$ be such that $t\gg H_1\gg 1$ and $H_2\gg 1$ .
At time $s=t/H_1$, we partition the space into two regions $D_1:=\{z:|z-x|>B_sH_2\}$ and $D_2:=\{z:|z-x|\leq B_sH_2\}$.
By the uniform bound and the large deviation,
\begin{equation}
	\P(\tau_x>t,X_s\in D_1,X_{t}=y)\leq \sum_{z\in D_1}\P(X^x_{s}=z)\P(X^z_{t-s}=y) \leq \frac{C}{B_tH_2^{\al-\eps}}.
\end{equation}
On $D_2$, we have $X^x_s=z$ with $|z|\leq |x|+B_sH_2$. 
Thus, by Propositions~\ref{P:dLLT}~and~\ref{P:dtLLT},
\begin{equation}\ba{l}
	\dis \P(\tau_x>t,X_s\in D_2,X^x_{t}=y) \leq \sum_{z\in D_2} \P_x(\tau_x>t,X_s=z)\P_z(X_{t-s}=y) \\ [5pt]
	\dis\qquad\leq \sum_{z\in \Z \setminus D_1}\P_x(\tau_x>t,X_s=z) \Big(\P_0(X_{t-s}=0)+\frac{C(|y|+|x|+B_sH_2)}{B_t^2}\Big) \\ [5pt]
	\dis\qquad\leq \Big(\P_x(\tau_x>t)+\frac{C}{H_2^{\al-\eps}}\Big) \Big(\P_0(X_{t}=0)+\frac{C}{B_tH_1}+\frac{C(|y|+|x|+B_sH_2)}{B_t^2}\Big)
\ea\end{equation}
Since $p(\cdot)$ has a finite $\kappa$-th moment for any $\kappa<1$, for any $\eps>0$,
\begin{equation}
	B_t^{-2}\sum_{|x|\leq B_t^2}|x|p(x)\leq B_t^{-2}\sum_{x} B_t^\eps\hs|x|^{1-\eps/2}p(x)\leq CB_t^{-2+\eps}.
\end{equation}
For the other error terms, fix any $\kappa<\al/(1-\al)$, and set $H_1=l(t)^{3/\al}\ell(t)^{\kappa}$ (so that $H_1\ll t$) and $H_2=H_1^{1/\al}$.
Then, using Lemma~\ref{L:jest}, direct calculations and basic estimates yield that
\begin{equation}
	\sum_xp(x)	\Big|\P_x(\tau_x>t,X_t=y)-\P_x(\tau_x>t)\P_0(X_t=0)\Big| \leq \frac{C }{B_tl(t)\ell(t)^{\kappa}}+\frac{C|y|}{B_t^2},
\end{equation}
where $\kappa$ can be made arbitrary close to $\al/(1-\al)$ as $\eps\to0$.
\epro

\section{First collision time of random walks}\label{S:fct}

Recall from Theorem~\ref{T:noncoa} that $(\Xb_t)_{t\geq0}=(X^1_t,\ldots,X^N_t)_{t\geq0}$ represents $N$ independent random walks with initial state $\Xb_0=\xb\in\Z^N$ and jump kernel $p(\cdot)$ in (\ref{jump}). 
Recall from (\ref{xtauc}) that $\tauc$ is the first collision time of $\Xb$.
We show Theorem~\ref{T:noncoa} for asymptotic probability of the non-collision event $\{\tauc>t\}$ in Subsection~\ref{S:noncoa}, and moreover show that this event is asymptotically independent of the transition event $\{X_t=\yb\}$ in Subsection~\ref{S:asyind}.

\subsection{Proof of Theorem~\ref{T:noncoa}}\label{S:noncoa}

As outlined in Subsection~\ref{S:ps}, our goal is to derive an ODE to characterize the asymptotic behavior of $\pc(t)=\P(\tauc>t)$.
The generator calculation yields (\ref{pc}), i.e.,
\begin{equation}\label{pcODE}
	-\frac{\di}{\di t}\pc(t)=\sum_{i<j}\sum_{y_{ij}\in\Z}2p(y_{ij}) \P(\tauc>t,\De^{ij}_t=y_{ij}),
\end{equation}
where $\De^{ij}_t=X^i_t-X^j_t$.
Our first step is to approximate the r.h.s.\ using the lemma below.
Recall the function class $\Bi$ such that $B_t\in\Bi$ implies $\ell(t)^K\geq l(t)$ for some $K>0$.
\bl \label{L:pcbdd1}
Assume that $B_t$ belongs to the class $\Bi$. Let $\kappa>0$.
Let $s=t/\ell(t)^{\theta}$ for some $\theta>K+\kappa+1$ if $\al=1$ or $\theta<\frac{1}{1-\al}$ if $\al<1$.
Then, there exist $c=c(\xb)$ and $\delta>0$ such that for any $\eps>0$,
\begin{equation}\label{tkl}
	\Big|\sum_{y_{ij}} p(y_{ij}) \P(\tau_\xb>t,\De^{ij}_t=y_{ij})-\frac{\pc(s)}{B_{2t}\hs l(2t)}\Big| \leq 
	\left\{\ba{ll}  
	\dis\frac{c }{B_{t}\hs l(t)} \Big(\frac{\pc(s)}{\ell(t)^{1-\eps}} +\frac{1}{\ell(t)^\kappa}\Big), & \al=1; \\ [10pt]
	\dis \frac{c}{B_{t}\hs l(t) t^\delta}, & \al\in(0,1).
	\ea\right.
\end{equation}
\el
\bpro
For any $k\neq l$, let $\tau_{kl}:=\inf\{t\geq0:X^k_t=X^l_t\}$ be the first collision time of $X^k$ and $X^l$. 
On $\{\tau_{\xb}>s,\De^{ij}_t=y_{ij}\}$, if $\{k,l\}\neq\{i,j\}$, then $\tau_{kl}\in(s,t]$ is unlikely.
Using the inclusion-exclusion principle, we obtain the following good approximation.
\begin{equation}\ba{l}\label{tkl1}
	\dis\quad\sum_{y_{ij}} p(y_{ij}) \big|\P(\tau_\xb>t,\De^{ij}_t=y_{ij})- \P(\tau_\xb>s,
	\tau_{ij}>t,\De^{ij}_t=y_{ij})\big| \\ [5pt]
	\dis \leq \sum_{\{k,l\}\neq\{i,j\}}\sum_{y_{ij}} p(y_{ij})\P(\tau_\xb>s,\tau_{ij}>t,\De^{ij}_t=y_{ij},s<\tau_{kl}\leq t).
\ea\end{equation}
We bound the r.h.s.\ of (\ref{tkl1}) by decomposing it into three cases as: {\bf (a)} $\{\tau_{kl}>s+u\}$, {\bf (b)}$\{\tau_{kl}\leq s+u, \hspace{2pt}|\De^{kl}_s|>B_u/\ell(u)^{\theta_1}\}$, {\bf (c)} $\{\tau_{kl}\leq s+u, \hspace{2pt} |\De^{kl}_s|\leq B_u/\ell(u)^{\theta_1}\}$. Here, we take $u=t/2$ if $\al=1$ and $u=t/\ell(t)^{\theta_2}$ otherwise, where $\theta_1>0,\theta_2\in(0,\frac{1}{1-\al})$ will be chosen later.

On the intersection with {\bf (a)}, we decompose at the collision time $\tau_{kl}.$
Let $v:=\tau_{kl}-s$.
For $\Xb_s=\zb\in\Z^N$, let $(\Xb_t^z)_{t\geq0}:=(\Xb_{t+s})_{t\geq0}$ denote the continuation of the random walks starting from time $s$.
We also use $\P_\zb$ to denote the law of random walks with a initial state $\zb$. 
Conditioning on $\Xb_s=\zb\in\Z^N$, the random walks first evolve to $\Xb^{\zb}_{v}=\textbf{w}$, with $\De_v^{kl}=w_{kl}=w_k-w_l$. 
Then at time $v$, $\De^{kl}$ jumps to 0 at rate $2p(w_{kl})$.
Over the remaining time $t-s-v$, $\De^{ij}$ moves from $w_{ij}$ to $y_{ij}$ without hitting 0.
Thus, the contribution from {\bf (a)} in the r.h.s.\ of (\ref{tkl1}) is bounded above by
\begin{equation}\ba{l}\label{tkl3}
	\dis \sum_{\{k,l\}\neq\{i,j\}}\int_{u}^{t-s} \sum_{\zb}\P(\tauc>s,\Xb_s=\zb) \sum_{\textbf{w}}\P_{\zb}(\tau_{kl}>v,\Xb_v=\textbf{w}) \\ [5pt] \dis\qquad\qquad\qquad\qquad\times2p(w_{kl})\sum_{y_{ij}}p(y_{ij})\P_{w_{ij}}(\tau_{ij}>t-s-v,\De^{ij}_{t-s-v}=y_{ij})\di v.
	\ea\end{equation}
Without loss of generality, assume that $i\neq k,l$. 
Let $X^{\{i\}}:=(X^1,\ldots,X^{i-1},X^{i+1},\ldots,X^N)$ denote the collection of random walks with the $i$-th coordinate omitted. We also denote $\textbf{w}^{\{i\}}:=(w^1,\ldots,w^{i-1},w^{i+1},\ldots,w^N)$.
We can then rewrite the sum over $\textbf{w}\in\Z^N $ as a sum over $\textbf{w}^{\{i\}}\in\Z^{N-1}$ and over $w_{ij}\in\Z$ (after replacing $w_i$ by $w_j+w_{ij}$).
Exploiting the independence between $X^i$ and other walks, and using the uniform bound for the transition probability,
\begin{equation}\label{tkl4}
	\P_\zb(\tau_{kl}>v,\Xb^{\{i\}}_v=\textbf{w}^{\{i\}},X^i_v=w_j+w_{ij})\leq \frac{C}{B_v} \P_\zb(\tau_{kl}>v,\Xb^{\{i\}}_v=\textbf{w}^{\{i\}}).
\end{equation}
Now for the event in the last sum in (\ref{tkl3}), it is equivalent to that the time-reversed random walk $\wih\De^{ij}$ moves from $y_{ij}$ to $w_{ij}$ during a time period $t-s-v=:r$ and $\hat\tau_{ij}>r$. Thus,
\begin{equation}\label{tkl5}
	\sum_{w_{ij}}\sum_{y_{ij}}p(y_{ij})\P_{w_{ij}}(\tau_{ij}>r,\De^{ij}_{r}=y_{ij})
	\dis=\sum_{y_{ij}}p(y_{ij})\P_{y_{ij}}(\hat\tau_{ij}>r) 	\leq \frac{C}{l(r)}.
\end{equation}
Inserting (\ref{tkl4})-(\ref{tkl5}) into (\ref{tkl3}) and using a similar time-reversal argument, we obtain an upper bound for the middle sum in (\ref{tkl3}) over $\textbf{w}^{\{i\}}$ as
\begin{equation}
	\sum_{\textbf{w}^{\{i\}}}p(w_{kl})\P_{\zb}(\tau_{kl}>v,\Xb^{\{i\}}_v=\textbf{w}^{\{i\}})=
	\sum_{w_{kl}}p(w_{kl})\P_{w_{kl}}(\hat\tau_{kl}>v,\wih\De^{kl}_{v}=z_{kl})
	\leq \frac{C}{B_vl(v)},
\end{equation}
where the inequality is due to Lemma~\ref{L:tau<}.
Consequently, (\ref{tkl3}) has an upper bound
\begin{equation}
	C_N \int_{u}^{t-s}\frac{1}{B^2_v\hs l(v)\hs l(t-s-v)} \sum_{\zb}\P(\tauc>s,\Xb=\zb) \di v
	= C_N \int_{u}^{t-s}\frac{\pc(s)}{B^2_v\hs l(v)\hs l(t-s-v)}\di v.
\end{equation}
Note that $t\gg s$ and $t/2\geq u$. Split the integral domain into $(u,2t/3)$ where $l(t-s-v)\asymp l(t)$, and $(2t/3,t-s)$ where $B_v^2\hs l(v)\asymp B_t^2\hs l(t)$. 
Applying the Potter bound and Karamata's theorem for the integrals then yields that for any $\eps>0$, (\ref{tkl3}) is bounded above by
\begin{equation}\label{tkl6}
	C_N\frac{t}{B^2_t\hs l(t)^2}\hs(t/u)^{\frac{2}{\al}-1+\eps}\pc(s)=\frac{C_N (t/u)^{\frac{2}{\al}-1+\eps}}{B_t\hs l(t)\ell(t)}\hs\pc(s).
\end{equation} 

On {\bf (b)}, during $(s+u,t]$, $\De^{ij}$ moves from some $w_{ij}\in\Z$ to $y_{ij}$ without hitting 0. By Markov's property the r.h.s.\ of (\ref{tkl1}) is
\begin{equation}\ba{l}
	\dis\sum_{\{k,l\}\neq\{i,j\}} \sum_{|z_{kl}|>B_u/\ell(u)^{\theta_1}}\sum_{w_{ij}}
	\P(\tau_\xb>s,s<\tau_{kl}\leq s+u,\De^{kl}_s=z_{kl}, \De^{ij}_{s+u}=w_{ij}) \\ [5pt]
	\hspace{3cm} \dis\times \sum_{y_{ij}} p(y_{ij})\P_{w_{ij}}(\tau_{ij}>t-s-u,\De^{ij}_{t-s-u}=y_{ij})
\ea\end{equation}
The time-reversal argument and Lemma~\ref{L:tauXt} imply that uniformly in $w_{ij}$,
\begin{equation}\label{tkl2}
	\dis\sum_{y_{ij}}p(y_{ij})\P_{w_{ij}}(\De^{ij}_{\ti r}=y_{ij},\tau>\ti r)
	=\dis\sum_{y_{ij}}p(y_{ij})\P_{y_{ij}}(\wih\De^{ij}_{\ti r}=w,\hat\tau>\ti r)
	\leq \frac{C}{B_t\hs l(t)},
\end{equation} 
where $\ti r:=t-s-u>t/3$.
Next, starting from $\De^{kl}_s=z_{kl}$ with $|z_{kl}|>B_u/\ell(u)^{\theta_1}$, the difference process $\De^{kl}$ must hit 0 during $(s,s+u]$. 
Hence, by Lemma~\ref{L:tau<}, for any $\eps>0$,
\begin{equation}
	\P_{z_{kl}}(\tau_{kl}\leq u)\leq C\ell(u)^{-1+(1-\al)\theta_1+\eps}.
\end{equation}
Applying the Markov property and using above bounds uniform in $w_{ij}$ and in $|z_{kl}|>B_u/\ell(u)^{\theta_1}$, we deduce that the contribution from {\bf (b)} in the r.h.s.\ of (\ref{tkl1}) is bounded above by 
\begin{equation}\label{tkl7}
	\sum_{\{k,l\}\neq\{i,j\}} \sum_{|z_{kl}|>B_u/\ell(u)^{\theta_1}} \P(\tau_\xb>s, \De_s^{kl}=z_{kl})\frac{C\ell(u)^{(1-\al)\theta_1+\eps}}{B_t\hs l(t)\ell(u)}
	\leq \frac{C_N\ell(u)^{(1-\al)\theta_1+\eps}}{B_t\hs l(t)\ell(u)}\pc(s).
\end{equation}

On {\bf (c)}, use the same time-reversal argument as for (\ref{tkl2}) for the interval $(s,t]$, and on $[0,s]$ simply apply the uniform transition probability for $\De^{kl}$. This yields the upper bound
\begin{equation}\ba{l}\label{tkl8}
	\dis\sum_{\{k,l\}\neq\{i,j\}} \sum_{z_{ij}}\P(|\De_s^{kl}|\leq B_u/\ell(u)^{\theta_1},\De_s^{ij}=z_{ij}) \sum_{y_{ij}}p(y_{ij})
	\P_{y_{ij}}(\wih\De^{ij}_{t-s}=z_{ij},\hat\tau_{ij}>t-s) \\ [5pt]
	\hspace{2cm}\dis\leq \frac{C_N B_u}{B_s\ell(u)^{\theta_1}}\hspace{1pt}\frac{1}{B_tl(t)}.
\ea\end{equation}

Combining the three bounds from (\ref{tkl6}), (\ref{tkl7}) and (\ref{tkl8}), we deduce that
\begin{equation}\label{tkl9}
	\text{r.h.s.\ of (\ref{tkl1})} \leq \frac{C_N\hs \pc(s)}{B_t\hs l(t)\ell(t)} \Big(\frac{\ell(t)}{\ell(u)^{1-(1-\al)\theta_1-\eps}} + (t/u)^{\frac{2}{\al}-1+\eps}\Big) + \frac{B_u}{B_s B_t\hs l(t)\ell(s)^{\theta_1}}.
\end{equation} 
For $\al=1$, recall $u=t/2$, and take $\theta_1=\kappa+2\theta$ so that $\frac{1}{B_s\ell(s)^{\theta_1}}\leq \frac{(t/s)^{1+\eps}}{B_t\ell(t)^{\theta_1}}\leq \frac{1}{B_t\ell(t)^{\kappa}}$ by the Potter bound. 
Inserting this into (\ref{tkl9}), basic estimates then imply the first line in the r.h.s.\ of (\ref{tkl}).
For $\al\in(0,1)$ and $u=t/\ell(t)^{\theta_2}$, recall that $\ell$ is $(1/\al-1)$-regularly varying. 
Bounding $\pc(s)\leq 1$ and taking suitable $\theta_1,\theta_2$, it is not hard to see that the second line in the r.h.s.\ of (\ref{tkl}) holds for some $\delta>0$.

To complete the proof, it remains to bound
\begin{equation*}\ba{l}
	\dis\quad \Big|\sum_{y_{ij}}p(y_{ij})\P(\tau_\xb>s,
	\tau_{ij}>t,\De^{ij}_t=y_{ij})-\frac{\pc(s)}{B_{2t}\hs l(2t)}\Big| \\ [5pt]
	\dis= \! \sum_{z_{ij}\in D_1}\!\!+\!\!\sum_{z_{ij}\in D_2}\!  \P(\tau_\xb>s,\De^{ij}_s=z_{ij})  \Big| \sum_{y_{ij}}p(y_{ij})\P_{y_{ij}}(\hat\tau_{ij}>t-s,\wih\De_{t-s}^{ij}=z_{ij}) - \frac{1}{B_{2t}\hs l(2t)}  \Big|,
\ea\end{equation*}
where $D_1:=\{z:|z|> B_s\ell(t)^{\theta_3}\}$ and $D_2=\{z:|z|\leq B_s\ell(t)^{\theta_3}\}$ with $\theta_3>0$ to be chosen later. To obtain the above, we have decomposed at time $s$, and considered the time-reversed rate 2 random walk $\wih\De^{ij}$ during $(s,t]$.
By using (\ref{tauxt}) and then large deviations, the sum for $z_{ij}$ over $D_1$ is bounded above by
\begin{equation}\label{bdd1}
	\frac{C}{B_t\hs l(t)}\hs\P(\De_s^{ij}\in D_2)\leq 
	\frac{Cs\hs L\big(B_s\ell(t)^{\theta_3}\big)}{B_t\hs l(t) B_s^\al\ell(t)^{\al\theta_3}}\leq \frac{C}{\Bl t\ell(t)^{\al\theta_3-\eps}}.
\end{equation}
On $D_2$, firstly, by (\ref{ptauXt}) in Lemma~\ref{L:tauXt} and by the fact that $|z_{ij}|\leq B_s\ell(t)^{\theta_3}$,
\begin{equation*}
	\sum_{y_{ij}}p(y_{ij}) \big| \P_{y_{ij}}(\hat\tau_{ij}>t-s,\wih\De_{t-s}^{ij}=z_{ij}) - \P_{y_{ij}}(\hat\tau_{ij}>t-s)\P_{0}(\wih\De_{t-s}^{ij}=0)\big|
	\leq \frac{C}{B_t\hs l(t)\ell(t)^{\theta_4}} + \frac{C\hs B_s\ell(t)^{\theta_3}}{B_t^2},
\end{equation*}
where we can choose any $\theta_4<\al/(1-\al)$.
Secondly, by Lemma~\ref{L:jest} and recalling from (\ref{bt}) the expression $\P_{0}(\wih\De_{t-s}^{ij}=0)=g_\al(0)/B_{2t-2s}$, we have
\begin{equation*}
	\dis \Big| \sum_{y_{ij}}p(y_{ij}) \P_{y_{ij}}(\hat\tau_{ij}>t-s)\P_{0}(\wih\De_{t-s}^{ij}=0) -\frac{1}{B_{2t}l(2t)} \Big|
	=\Big| \frac{1+O(1/\ell(t))}{g_\al(0)\hs l(2t-2s)}\frac{g_\al(0)}{B_{2t-2s}}-\frac{1}{B_{2t} l(2t)} \Big|.
\end{equation*}
We then insert $\frac{1}{l(2t-2s)}=\frac{1}{l(2t)}+O(\frac{1}{l(t)\ell(t)})$ by (\ref{dl}) and insert $\frac{1}{B_{2t-2s}}=\frac{1}{B_{2t}}+O(\frac{s}{tB_t})$ by Proposition~\ref{P:dtLLT} in the above equality.
This, together with the error bound in the first step, yields that the sum for $z_{ij}$ over $D_2$ is bounded above by
\begin{equation}\label{bdd2}
C\hs\pc(s)\Big( \frac{1}{B_t\hs l(t)\ell(t)^{\theta_4}} + \frac{B_s\ell(t)^{\theta_3}}{B_t^2} +\frac{1}{B_t\hs l(t)\ell(t)} +\frac{(s/t)}{B_t\hs l(t)} \Big)
\end{equation}
with $\theta_4<\al/(1-\al)$.
For $\al=1$, choose, e.g., $\theta_3=\kappa+1/2$ and $\theta_4=1$.
Since $B_t$ is in the class $\Bi$, for $s=t/\ell(t)^{\theta}$ with $\theta>K+\kappa+1$, by basic estimates the bounds (\ref{bdd1}) and $(\ref{bdd2})$ can be simplified to
\begin{equation*}
	\frac{C}{B_t\hs l(t)}\Big(\frac{\pc(s)}{\ell(t)^{1-\eps}} + \frac{1}{\ell(t)^\kappa}\Big).
\end{equation*}
For $\al<1$, recall that $l$ is slowly varying and $\ell$ is $(1/\al-1)$-regularly varying.
Then, for any $0<\theta_3<\theta/\al$ and $\theta_4>0$, there exists some $\de>0$ such that (\ref{bdd1}) and $(\ref{bdd2})$ can be upper bounded by 
\begin{equation*}
	\frac{C}{B_t\hs l(t)t^{\de}}.
\end{equation*}
The proof is thus finished.
\epro

\bl\label{L:pcbdd}
Recall that $\ups_N={N\choose 2}$. Under the assumptions of Theorem~\ref{T:noncoa}, there exists $C_\xb$ such that
\begin{equation}\label{pcbdd}
	\P(\tauc>t)\leq \frac{C_\xb}{l(t)^{\ups_N}} \quand
	\Big|\frac{\di}{\di t}\P(\tauc>t)\Big|\leq \frac{C_\xb}{B_t\hs l(t)^{\ups_N+1}}.
\end{equation}
\el
\bpro
For $\al=1$, note that $l(t)\leq \ell(t)^{K}$ as $B_t\in\Bi$.
We take $\kappa=K(\ups_N+1)$ in Lemma~\ref{L:pcbdd1} and then insert (\ref{tkl}) into the ODE (\ref{pcODE}).
This yields that for any $\eps>0$,
\begin{equation}\label{pcODE2}
	-\frac{\di}{\di t}\pc(t) 
	\geq\frac{2\ups_N - C\hs\ell(t)^{-1+\eps}}{B_{2t}l(2t)}\pc(t) -\frac{C}{B_tl(t)^{\ups_N+2}},
\end{equation}
where we used $\pc(s)\geq \pc(t)$. Define $F(t):=\int_{0}^{t}\frac{C}{B_{2s}\hs l(2s)\ell(s)^{1-\eps}}\di s$, which is bounded since 
\begin{equation*}
	\int_0^t\frac{\di s}{B_{2s}\hs l(2s)\ell(s)^{1-\eps}}\leq \int_0^\infty\frac{C\hs\di l(s)}{\hs l(s)^{1+(1-\eps)/K}}<\infty.
\end{equation*}
As $\big(\log(l(2t)^{\ups_N})\big)'=\frac{2\ups_N}{B_{2t}l(2t)}$, multiplying to (\ref{pcODE2}) the integrating factor $\exp\big(\log(l(2t)^{\ups_N})-F(t)\big)=e^{-F(t)}l(2t)^{\ups_N}$ and using $e^{-F(t)}\leq 1$, we obtain
\begin{equation*}\label{dpcbdd}
	\big(e^{-F(t)} l(2t)^{\ups_N} \pc(t)\big)' 
	\leq \frac{C}{B_t\hs l(t)^{2}}.
\end{equation*}
The r.h.s.\ has a bounded integral. This and boundedness of $F(t)$ imply the first inequality in (\ref{pcbdd}).

To prove the second inequality in (\ref{pcbdd}), we need to derive a bound in the opposite direction to that of~(\ref{pcbdd}) for the derivative of $\pc(t)$. 
Similar calculations as in the proof of Lemma~\ref{L:dj} yield
\begin{equation}\label{lNd}
	|l(t)^{-\ups_N}-l(s)^{-\ups_N}|\leq \frac{C_N (t/s)^\eps}{l(t)^{\ups_N}\ell(t)}.
\end{equation}
Now for $s=t/\ell(t)^{\theta}$ in Lemma~\ref{L:pcbdd1}, using the above bound and $\pc(s)\leq C\hs l(2s)^{-\ups_N}$, and applying (\ref{tkl}) into the ODE (\ref{pcODE}), we obtain
\begin{equation}
	\Big|\frac{\di}{\di t}\P(\tauc>t)\Big|\leq \frac{C_\xb}{B_t\hs l(t)^{\ups_N+1}}.
\end{equation}

The same method with simpler calculations is applicable for $\al\in(0,1)$, and this completes the proof.
\epro

\bpro[of Theorem~\ref{T:noncoa}]
Take $s=t/\ell(t)^{\theta}$ as in Lemma~\ref{L:pcbdd1}. Then, by Lemma~\ref{L:pcbdd} and (\ref{lNd}),
\begin{equation}\label{taust}
	\P(\tauc\in(s,t])\leq \int_s^t\frac{C_\xb}{B_u\hs l(u)^{\ups_N+1}}\di u
	\leq \frac{C_\xb}{l(t)^{\ups_N}\ell(t)^{1-\eps}}.
\end{equation}
Therefore, by (\ref{pcODE}) and (\ref{tkl}), there exists some $\de>0$ such that
\begin{equation*}
	-\frac{\di}{\di t}\P(\tauc>t)=\frac{2\ups_N}{B_{2t}l(2t)}\P(\tauc>t) + R(t),
\end{equation*}
where the error term $R(t)=O\Big(\frac{C_\xb}{B_t\hs l(t)^{\ups_N+1}\ell(t)^{1-\eps}}\Big)$ if $\al=1$ and $R(t)=O\Big(\frac{C_\xb}{B_t\hs t^{\de}}\Big)$ if $\al<1$.

For $\al=1$, since $\frac{1}{t\hs\ell(t)^{2-\eps}}=\frac{1}{B_t\hs l(t)\ell(t)^{1-\eps}}\leq \frac{1}{B_t\hs l(t)^{1+\eps/K}}$ is integrable, by Karamata's theorem $\int_{t}^{\infty}\frac{\di u}{u\hs\ell(u)^{2-\eps}}\to0$ is slowly varying and $\frac{1}{\ell(t)^{2-\eps}}\int_{t}^{\infty}\frac{\di u}{u\hs\ell(u)^{2-\eps}}\to\infty$.
We then have
\begin{equation*}
	\P_\xb(\tauc>t)\cdot l(2t)^{\ups_N}=\int_{0}^\infty h(u)\di u -\int_t^\infty h(u)\di u,
\end{equation*}
where $h(t)=l(2t)^{\ups_N}R(t)=O\Big(\frac{C_\xb}{t\hs \ell(t)^{2-\eps}}\Big)$.
Hence, setting $c_\xb=\int_{0}^\infty h(u)\di u$, we have
\begin{equation}\label{pnc1}
	\P_\xb(\tauc>t)=\frac{c_\xb}{l(2t)^{\ups_N}}\Big(1+O(\Ri_{\rm NC}(\xb,t))\Big),
\end{equation}
where $\Ri_{\rm NC}(\xb,t)=C_\xb\int_{t}^{\infty}\frac{\di u}{u\hs \ell(u)^{2-\eps}}= o\big(\ell(t)^{-2+\eps}\big)$. 

For $\al<1$, recall that $B_t$ is $1/\al$-regularly varying. 
The same argument implies that (\ref{pnc1}) holds with $\Ri(t)=C_\xb t^{-\eta}$, where $\eta$ is a constant in $(1/\al-1,\hs \de+1/\al-1)$.
\epro

\subsection{Asymptotic independence}\label{S:asyind}
For $N$ independent random walks $(\Xb_t)_{t\geq0}$, recall that we have also adopted the convention that $\P_\xb(\cdot)$ denotes the law for $\Xb$ with initial state $\Xb_0=\xb$.
In this subsection, we aim to prove the asymptotic independence between the transition $\{\Xb_t=\yb\}$ of the random walks and the non-collision event $\{\tauc>t\}$ in the following sense.
\begin{proposition}[Asymptotic independence]\label{L:retran}
Under the assumptions of Theorem~\ref{T:noncoa}, for independent random walks $(\Xb_t)_{t\geq0}$, there exist $C_\xb$ and $\de>0$ such that for any $\eps>0$ and sufficiently large $t$,
\begin{equation}\label{retran1}
	\sum_{\yb\in\Z^N}\Big| \P_\xb(\Xb_t=\yb,\tauc>t)-\P_\xb(\Xb_t=\yb) \P_\xb(\tauc>t)\Big|\leq 
	\left\{\ba{ll}  
	\dis\frac{C_\xb}{l(t)^{\ups_N}\hs\ell(t)^{1-\eps}}, & \al=1; \\ [10pt]
	\dis \frac{C_\xb}{l(t)^{\ups_N}\hs t^{\de}}, & \al\in(0,1).
	\ea\right.
\end{equation}
\end{proposition}
\bpro
Let $s=t/\ell(t)^{\theta_1}$ with $\theta_1<\frac{\al}{1-\al}$. 
For the l.h.s.\ of (\ref{retran1}), we bound the summand for each $\yb$ by the sum of the following four differences.
\begin{equation*}\ba{l}
	\dis R_1^\yb := \big|\P_\xb(\Xb_t=\yb,\tauc>t)-\P_\xb(\Xb_t=\yb,\tauc>s)\big| , \\ [5pt]
	\dis R_2^\yb := \big|\P_\xb(\Xb_t=\yb) \P_\xb(\tauc>t)-\P_\xb(\Xb_t=\yb) \P_\xb(\tauc>s)\big|, \\ [5pt]
	\dis R_3^\yb := \dis \big|\P_\xb(\Xb_t=\yb,\tauc>s)-\P_\xb(\Xb_{t-s}=\yb)\P_\xb(\tauc>s)\big|, \\ [5pt]
	\dis R_4^\yb := \dis \big|\P_\xb(\Xb_t=\yb)\P_\xb(\tauc>s)-\P_\xb(\Xb_{t-s}=\yb)\P_\xb(\tauc>s)\big|.
\ea\end{equation*}
\noi{\bf Bounding $\sum_{\yb\in\Z^N} (R^\yb_1+R^\yb_2)$.} Since
\begin{equation*}
	R_1^\yb = \P_\xb(\Xb_t=\yb,\tauc\in(s,t]) \quand  
	R_2^\yb = \P_\xb(\Xb_t=\yb)\P_\xb(\tauc\in(s,t]),
\end{equation*}
The same calculations as in (\ref{taust}) lead to the fact that for any $\eps>0$,
\begin{equation}\label{R12}
	\sum_{\yb\in\Z^N}(R^\yb_1+R^\yb_2)=2\P_\xb(\tauc\in(s,t])\leq  \frac{C_\xb}{l(t)^{\ups_N}\ell(t)^{1-\eps}}.
\end{equation}
\noi {\bf Bounding $\sum_{\yb\in\Z^N} (R^\yb_3+R^\yb_4)$.}
First, we decompose at time $s$ and apply the Markov property,
\begin{equation}\label{R31}
	\sum_{\yb} R_3^\yb\leq \sum_{\yb}\sum_\zb\P_\xb(\Xb_s=\zb,\tauc>s)\big|\P_{\zb}(\Xb_{t-s}=\yb)-\P_{\xb}(\Xb_{t-s}=\yb)\big|.
\end{equation}
A similar bound also holds for $R_4^\yb$. We then bound the r.h.s.\ of (\ref{R31}).

For $\xb=(x_1,\ldots,x_N)$, denote $|\xb|_{\rm m}:=\max_{1\leq i\leq N}|x_i|$.
Let $\theta_2>0$. 
Let $D_1:=\{\wb\in\Z^d:|\wb-\xb|_\mr\leq B_t\hs \ell(t)^{\theta_2}\}$ and $D_2:=\{\wb\in\Z^d:|\wb-\xb|_\mr\leq B_s\hs \ell(t)^{\theta_2}\}$.
We consider three regions: {\bf (a)} $\{\Xb_t\in D_1^c\}$,  {\bf (b)} $\{\Xb_s\in D_2^c\}$, and {\bf (c)} $\{\Xb_t\in D_1,\Xb_s\in D_2\}$.
The large deviation (\ref{ld}) and basic estimates, such as the Potter bound and (\ref{btllt}), give a bound of the probabilities of the first two events {\bf (a)} and {\bf (b)}:
\begin{equation}\label{R32}
	\P(\Xb_t\in D_1^c) +\P(\Xb_s\in D_2^c)\leq \frac{C_\xb}{\ell(t)^{\al\theta_2-\eps}}.
\end{equation}
On {\bf (c)}, using the union bound $C/B_t$ for the transition probability, and applying the inequality $|\prod_{i=1}^{N} a_i-\prod_{i=1}^{N} b_i|\leq \sum_{i=1}^{N}\prod_{j<i}a_j \hs |a_i-b_i| \prod_{j>i}b_j\leq C^{N-1}\sum_{i=1}^{N}|a_i-b_i|$ for $0\leq a_i,b_i\leq C$,
\begin{equation*}
	\big|\P_{\zb}(\Xb_{t-s}=\yb)-\P_{\xb}(\Xb_{t-s}=\yb)\big|\leq (C/B_t)^{N-1} \sum_{i=1}^{N}\big|\P_{z_i}(X^i_{t-s}=y_i)-\P_{x_i}(X^i_{t-s}=y_i)\big|.
\end{equation*}
Therefore, by Proposition~\ref{P:dLLT} and $|z_i-x_i|\leq B_s\hs\ell(t)^{\theta_2}$ for $\zb\in D_2$, we have
\begin{equation}\ba{l}\label{R33}
	\dis\sum_{\yb\in D_1}\sum_{\zb\in D_2}\P_\xb(\Xb_s=\zb,\tauc>s) \big|\P_{\zb}(\Xb_{t-s}=\yb)-\P_{\xb}(\Xb_{t-s}=\yb)\big| \\ [5pt]
	\dis\quad\leq \frac{C_\xb  B_s\hs\ell(t)^{\theta_2}}{B_t^{N+1}}\sum_{\yb\in D_1}\sum_{\zb\in D_2}\P_\xb(\Xb_s=\zb,\tauc>s)
	\leq \frac{C_\xb  B_sB_t^N\hs\ell(t)^{(N+1)\theta_2}}{B_t^{N+1}\hs l(t)^{\ups_N}},
\ea\end{equation}
where in the last inequality we used that the sum over $\zb$ is bounded by $\P(\tau_\xb>s)\leq C_\xb l(t)^{-\ups_N}$ thanks to Lemma~\ref{L:pcbdd}.
For $\al=1$, if we take $\theta_2=K\ups_N+1$ and $\theta_1=(N+1)\theta_2+1$, then combining the three bounds (\ref{R12}), (\ref{R32}) and (\ref{R33}) we arrive at the desired bound (\ref{retran1}).
For $\al<1$, we take $\theta_2$ sufficiently small and $\theta_1=(N+1)\al\theta_2+\al$ such that $\theta_1<\frac{\al}{1-\al}$ is satisfied. Then, (\ref{R12}), (\ref{R32}) and (\ref{R33}) imply that there exists some $\de$ such that (\ref{retran1}) holds, which finishes the proof.
\epro
\br\label{R:asyind}
Via a simple adaption of the proof above, we can change the fixed $\xb$ to the origin ${\bf 0}$ such that the same upper bound as in (\ref{retran1}) holds for
\begin{equation*}
	\sum_{\yb\in\Z^N}\Big| \P_\xb(\Xb_t=\yb,\tauc>t)-\P_{\bf 0}(\Xb_t=\yb) \P_\xb(\tauc>t)\Big|.
\end{equation*}
Indeed, comparing with (\ref{ptauXt}) in Lemma~\ref{L:tauXt}, we can replace the term $\P_\xb(\Xb_{t-s}=\yb)$ in (\ref{R31}) by $\P_{\bf 0}(\Xb_{t-s}=\yb)$, and obtain the same upper bound as in (\ref{R33}).
Thus, the desired result follows once $\P_\xb(\Xb_t=\yb)$ and $\P_\xb(\Xb_{t-s}=\yb)$ in $R^\yb_2$ and $R^\yb_4$ are replaced by $\P_{\bf 0}(\Xb_t=\yb)$ and $\P_{\bf 0}(\Xb_{t-s}=\yb)$, respectively.
\er

\section{Coalescing heavy-tailed random walks}\label{S:npt}

Our goal is to show Theorem~\ref{T:den} for the density function asymptotics in Subsection~\ref{S:thmden}, and show Theorem~\ref{T:npf} for the $N$-point correlation function asymptotics in Subsection~\ref{S:thmnpt}.
As a preliminary, we give an estimate of the $N$-point function in Subsection~\ref{S:nptest}.

\subsection{$N$-point correlation function}\label{S:nptest}

For $\xb\in\Z^N$, recall the $N$-point function $\rho_N(\xb,t)$ from (\ref{npt2}) and recall the density function $\rho_1(t)$.
For the CRW $(\xi_t)_{t\geq0}$ starting from every site on $\Z$ and for any $\wb=(w_1,\ldots,w_n)\in\Z^n$, with a bit abuse of notation, we define $\xi_t(\wb):=\prod_{i=1}^{n}\xi_t(w_i)\in\{0,1\}$.
Then, we can rewrite that $\rho_N(\xb,t)=\E[\xi_t(\xb)]$.
In this subsection, we denote $\pc(t):=\P(\hat\tau_{\bf x}>t)$.
\begin{lemma}\label{L:nptest}
Let $s< t$. Under the assumptions of Theorem~\ref{T:npf}, there exists $\de>0$ such that
\begin{equation}\label{Ri}
	\E[\xi_t(\xb)]=\rho_1(t-s)^N\hs \pc(s)\big( 1+\Ri_1+\Ri_2+\Ri_3 \big),
\end{equation}
where $0\geq\Ri_1=O\Big(\sum\limits_{J=1}^N\big(\frac{s\hs \rho_1(t-s)}{l(s)}\big)^J\Big)$, $\Ri_2=O\Big( \frac{1}{\ell(t)^{1-\eps}\wedge t^\de} \Big)$ and $0\leq\Ri_3=O\Big( \sum\limits_{J=1}^{N-1}\big(\frac{1}{B_s\hs \rho_1(t-s)}\big)^J \Big)$.
\end{lemma}
\br
Later when we use Lemma~\ref{L:nptest}, we choose $s$ such that $\frac{s\hs \rho_1(t-s)}{l(s)},\frac{(t/s)^\eps}{B_t\hs \rho_1(t-s)}\ll1$, and thus $\Ri_1$ and $\Ri_3$ can be written in much tighter forms.
\er
\bpro
Firstly, we trace back to time $t-s$ for the ancestors of the $N$ particles at $\xb$ at time $t$.
For $\yb\in\Z^N$, we say that $(\yb,t-s)$ connects to $(\xb,t)$ if $N$ coalescing random walks in the system move from $\yb$ at time $t-s$ to $\xb$ at $t$.
By the natural coupling between the coalescing random walks and independent random walks, and by the time-reversal argument,
\begin{equation*}
	\E[\xi_t(\xb)]=\sum_{\yb\in\Z^N}\E[\xi_{t-s}(\yb)]\P_\xb(\wih\Xb_s=\yb,\hat\tau_{\rm \xb}>s)+R_1,
\end{equation*}
where $(\wih\Xb_t)_{t\geq0}$ denotes $N$ time-reversed independent random walks starting from $\xb$, $\hat\tau_{\rm \xb}$ is their first collision time, and $R_1$ is the error term. 
We remark that if $y_i=y_j$ for some $i\neq j$ in $\yb$, then $\P_\xb(\wih\Xb_s=\yb,\hat\tau_{\rm \xb}>s)=0$.

Indeed, if we call $A(\yb)$ the event that in the coalescing system there are particles at $\yb$ at time $t-s$ which connects to $\xb$ at time $t$, then in the equation above, the l.h.s.\ equal $\P(\cup_\yb A(\yb))$ and the sum in the r.h.s.\ equals $\sum_\yb\P(A(\yb))$.
Thus by the inclusion-exclusion principle, $R_1\leq0$ and $|R_1|$ is bounded above by the probability of the following event: there exist $\yb\neq\zb\in\Z^N$ such that {\bf (a)}~$\xi_{t-s}(\yb,\zb)=1$ (i.e., there are particles at all sites of $\{\yb,\zb\}$ at time $t-s$), and {\bf (b)} both $(\yb,t-s)$ and $(\zb,t-s)$ connect to $(\xb,t)$.
Suppose that there are $J$ different coordinates of $\yb$ and $\zb$ and the remaining $N-J$ coordinates are the same.
W.l.o.g., we rearrange the order of coordinates, and assume that the first $J$ coordinates of $\yb$ and $\zb$ differ.
Since then the two random walks starting from $y_i$ and $z_i$ at $t-s$ must coalesce before time $t$ for $1\leq i\leq J$, assume additionally that their coalescing times are ordered $t-s<u_J<\ldots<u_1<t$.
In words, {\bf (a)} is the event that in the coalescing system there are $N+J$ particles occupying the sites $\{\yb,\zb\}=\{y_1,\ldots,y_N,z_1,\ldots,z_J\}$ at time $t-s$, and {\bf (b)} is the event that the pair of particles at site $y_i$ and $z_i$ move to $x_i$ at time $t$ for $1\leq i\leq J$ and the particle at $y_j$ moves to $x_j$ at time $t$ for $J<j\leq N$.
$|R_1|$ is then bounded from above by summing the probabilities of {\bf (a)} intersects {\bf (b)} for all $\{\yb,\zb\}$ over $\Z^{N+J}$ and then for $J$ from $1$ to $N$.

To obtain an upper bound for the probability of {\bf (a)}, we recall that there are exactly $N+J$ particle at $\{\yb,\zb\}$ at $t-s$, and by the negative correlation of CRW (see Theorem~\ref{T:necor}),
\begin{equation}\label{Ri12}
	\E[\xi(\yb,\zb)]\leq \rho_1(t-s)^{N+J}.
\end{equation}
We then take {\bf (b)} into account. 
It implies that the following event on the time-reversed (branching) random walks occurs.
Set $\hat u_i:=t-u_i$.
For the reversed time from 0 to $\hat u_1$, $N$ random walks $(\wih\Xb_t)_{t\geq0}$ starting from $\xb$ does not collide.
For $1\leq i\leq J$, from $\hat u_i$ to $\hat u_{i+1}$ (with $\hat u_{J+1}:=s$), we compare with the forward coalescing system. 
Given $\wih X^i_{\hat u_i}$, there is a binary branching at time $\hat u_i$ due to the coalescence in the forward coalescing system, where one descendant stays and the other jumps to $\wih Y^i_{\hat u_i}=\wih X^i_{\hat u_i}+w_i$ with probability $p(w_i)$, and the pair does not collide at least by time $\hat u_{i+1}$.
We note that the branching rate is $2$ as in the forward system either particle in $i$-th pair can jump and coalesce with the other.
Since (\ref{Ri12}) gives a uniform bound, for the upper bound of $|R_1|$, the sum of $\{yb,\zb\}$ is applied to {\bf (b)} only, and it is equivalent to that there is no restriction for the time-reversed random walks at reversed time $s$.
This yields an upper bound of $|R_1|$ as
\begin{equation*}
	\sum_{J=1}^{N}\rho_1(t-s)^{N+J} \hspace{-15pt} \idotsint\displaylimits_{0<\hat u_1<\hat u_2<\ldots<\hat u_J<s} \hspace{-5pt} C_\xb \hspace{-5pt}\sum_{w_1,\ldots,w_J}\hspace{-5pt}\P_\xb(\htauc>\hat u_1)p(y_1) \P_{y_1}(\tau_{y_1}>\hat u_1)\ldots p(y_J) \P_{y_J}(\tau_{y_J}>\hat u_J)\di\hat u_1\ldots\di\hat u_J.
\end{equation*}
By Theorem~\ref{T:noncoa} and Lemma~\ref{L:jest}, the multiple integral in the above display is bounded from above by
\begin{equation}
	C_\xb\hs \P_\xb(\htauc>s) \big(s/l(s)\big)^J = C_\xb\hs \pc(s) \hs\big(s/l(s)\big)^J.
\end{equation}
Therefore, we have
\begin{equation*}
	|R_1| \leq C_\xb\sum_{J=1}^N \rho_1(t-s)^{N+J} \pc(s) \hs\Big(\frac{s}{l(s)}\Big)^J 
	= \rho_1(t-s)^N\hs\pc(s) O\Big( \sum_{J=1}^{N}\Big(\frac{s\hs \rho_1(t-s)}{l(s)}\Big)^J\Big),
\end{equation*}
which gives rise to $\Ri_1$ in (\ref{Ri}).

Secondly, we write
\begin{equation*}
	\sum_{\yb\in\Z^N}\E[\xi_{t-s}(\yb)]\P_\xb(\wih\Xb_s=\yb,\hat\tau_{\rm \xb}>s)
	= \sum_{\yb\in\Z^N}\E[\xi_{t-s}(\yb)]\P_\0(\wih\Xb_s=\yb)\P_\xb(\htauc>s)+R_2.
\end{equation*}
By the negative correlation $\E[\xi_{t-s}(\yb)]\leq \rho_1(t-s)^{N}$, the asymptotic independence Proposition~\ref{L:retran} and Remark~\ref{R:asyind} then imply that there exists $\de>0$ such that
\begin{equation}
	|R_2|\leq \frac{C_\xb \hs \rho_1(t-s)^N}{l(s)^{\ups_N}(\ell(s)^{1-\eps}\wedge t^{\de})}
	\leq \rho_1(t-s)^N\hs\pc(s) \frac{C_\xb}{(\ell(s)^{1-\eps}\wedge t^{\de})},
\end{equation}
where we used Theorem~\ref{T:noncoa} for $\pc$ in the last inequality.
This leads to $\Ri_2$ in (\ref{Ri}).

Thirdly, we approximate 
\begin{equation*}
	\sum_{\yb\in\Z^N}\E[\xi_{t-s}(\yb)]\P_\0(\wih\Xb_s=\yb)\P_\xb(\htauc>s)
	=\prod_{i=1}^N \Big(\sum_{y_i\in\Z}\E[\xi_{t-s}(y_i)]\P_0(\wih X^i_s=y_i) \Big)\P_\xb(\htauc>s) +R_3.
\end{equation*}
Denote $Z:=\sum_i\xi_{t-s}(y_i)\P_0(\wih X_s^i=y_i)$.
By Jensen's inequality $\E[Z^N]\geq (\E[Z])^N$, we see that $R_3\geq0$.
For the other inequality, let $\Z^N_J:=\{\yb\in\Z^N:J=|\{y_1,\ldots,y_N\}|\}$ be the set of those $\yb$'s with $J$ different coordinate values.
For $J=N$, by the negative correlation we have
\begin{equation*}
	\sum\nolimits_{\yb\in\Z^N_N}\big(\E[\xi_{t-s}(\yb)]-\prod\nolimits_{i}\E[\xi_{t-s}(y_i)]\big)\P_{\bf 0}(\wih\Xb_s=\yb) \leq 0
\end{equation*}
For $1\leq J<N$, w.l.o.g.\ we assume that $y_1,\ldots,y_J$ are different and are repeated in $\yb$ for $k_1,\ldots,k_J$ many times, respectively.
Then, by the negative correlation and by the uniform bound for the transition probability,
\begin{equation*}
	\sum_{\yb\in\Z^N_J}\E[\xi_{t-s}(\yb)] \P_{\bf 0}(\wih\Xb_s=\yb)\leq C\!\!\!\sum_{y_1,\ldots,y_J} \rho_1(t-s)^{J}\prod_{i=1}^{J} \frac{\P_0(\wih X_s^i=y_i)}{B_s^{k_i-1}}=\frac{C\hs \rho_1(t-s)^J}{B_s^{N-J}},
\end{equation*}
where in the equality we used $\sum_{y_i}\P_0(\wih X_s^i=y_i)=1$ and $\sum_i k_i=N$.
Therefore,
\begin{equation*}
	R_3\leq C_\xb\hs \rho_1(t-s)^{N}\hs \pc(s) \sum_{J=1}^{N-1} \frac{1}{(\rho_1(t-s)B_s)^J},
\end{equation*}
which gives $\Ri_3$.
The proof is thus complete. 
\epro

\subsection{Proof of Theorem~\ref{T:den}}\label{S:thmden}

To give precise asymptotics of the density, we need a priori rough bounds of $\rho_1$ that give the correct order, where the upper bound can be proved using a modification of the arguments of \cite{BG80a}.

\begin{proposition}\label{T:priori}
	Under the assumptions of Theorem~\ref{T:npf}, there exists $c,C>0$ such that for all large $t$,
	\begin{equation}\label{pribdd2}
		\frac{c\hs l(t)}{t}\leq\rho_1(t)\leq\frac{C\hs l(t)}{t}.
	\end{equation}
\end{proposition}
\bpro
The proof proceeds in three steps, establishing successively a crude lower bound $\rho_1(t)\ge c/t$, an upper bound $\rho_1(t)\le C\hs l(t)/t$, and finally an improved lower bound $\rho_1(t)\ge c\hs l(t)/t$.

{\it Step 1.} Recall the notation that $\xi(0,x)=\xi(0)\xi(x)$.
By the generator calculation, using the negative correlation $\E[\xi(0,x)]\leq\E[\xi(0)]\E[\xi(x)]$ (see Theorem~\ref{T:necor}), and by the translation invariance $\E[\xi(x)]=\rho_1(t)$,
\begin{equation}\label{pribdd}
	-\frac{\di}{\di t}\rho_1(t)=\sum_{x\in\Z}p(-x)\E[\xi_t(0)\xi_t(x)]\leq \sum_{x\in\Z}p(-x)\rho_1(t)^2=\rho_1(t)^2,
\end{equation}
which yields a crude bound
\begin{equation*}
	\rho_1(t)\geq \frac{c}{t}.
\end{equation*}

{\it Step 2.} The idea of proving the upper bound is credited to \cite{BG80a}.
We thus only provide a proof sketch here.
The idea is that given the density $\rho_1(t)$, starting from the interval of length $C/\rho_1(t)$, there is a probability bounded uniformly from below that all other particles would coalesce with one particle over time $s$, which results to a lower bound for the density decay.
For any $L\in\N$, Denote $\Lambda:=[0,L]\cap\Z$ and $x+\Lambda:=\{x+y:y\in\Lambda\}$. Let $\xi^A_t$ be the coalescing random walks starting from set $A$. Identifying $\xi_t(B)$ as the set $\{x\in B:\xi_t(x)=1\}$ for any set $B\subseteq\Z$, we have that for any $s,t\geq0$,
\begin{equation*}
	\E\big[|\xi_{t+s}(\Lambda)|\big]
	\leq \sum_{x\in L\Z}\E\big[|\xi^{\xi_t(x+\Lambda)}_s(\Lambda)|\big]
	=\sum_{x\in L\Z}\E\big[|\xi^{\xi_t(\Lambda)}_s(x+\Lambda)|\big]
	=\E\big[|\xi^{\xi_t(\Lambda)}_s(\Z)|\big],
\end{equation*}
where in the inequality we applied the Markov property at time $t$ for ancestors of the particles in $\Lambda$ at time $t+s$, and in the first equality we used translation invariance.
We note that all terms in the equation above are finite as $\E\big[|\xi^{\xi_t(\Lambda)}_s(\Z)|\big]\leq \E\big[|\xi_t(\Lambda)|\big]\leq |\Lambda|$.
Since $\E[|\xi_t(\Lambda)|]=L\rho_1(t)$, we derive 
\begin{equation}\label{rhobdd}
	L(\rho_1(t)-\rho_1(t+s))=\E[\xi_t(\Lambda)]-\E\big[|\xi_{t+s}(\Lambda)|\big]\geq \E[\xi_t(\Lambda)-\xi^{\xi_t(\Lambda)}_s(\Z)\big].
\end{equation}
The right-hand side of (\ref{rhobdd}) is the expected number of times that the particles $\xi_t(\Lambda)$ coalesce through time $s$.
This is lower bounded if we fix a particle in $\xi_t(\Lambda)$, and consider the number of coalescence between this particle and others only, which is further bounded from below by
\begin{eqnarray*}\label{rhobdd2}
	\E\big[\xi_t(\Lambda)-1\big]\cdot\min_{|x|\leq L}\P(\tau_x\leq 2s).
\end{eqnarray*}
Recall that $\De(x,s)=\P(X^0_s=0)-\P(X^x_s=0)$.
Using the last-exit decomposition (\ref{xexit}),
\begin{equation*}
	\P(\tau_x\leq s)= 1- \Delta(x,s) - \int_0^s \De(x,u)j(s-u)\di u.
\end{equation*}
Since $|\De(x,u)|\leq\frac{c|x|}{B_u^2}\leq\frac{cL}{B_u^2}$ and $j(u)\sim 1/l(u)$, 
it follows $|\int_0^s \De(x,u)j(s-u)\di u|\leq \frac{csL}{B_s^2l(s)}$.
Inserting all these bounds into (\ref{rhobdd}), we get
\begin{equation}\label{drho}
	L\rho_1(t)-L\rho_1(t+s)\geq \big(L\rho_1(t)-1\big)\Big(1-\frac{csL}{B_s^2l(s)}\Big),
\end{equation}
implying
$L\rho_1(t+s)\leq P L \rho_1(t) +1,$
where $P:=csL/\big(B_s^2l(s)\big)$.
Repeating this for $\kappa:=\lfloor t/s\rfloor+1$ many times yields
\begin{equation*}
	L\rho_1(2t)\leq P^\kappa L\rho_1(t)+(1-P)^{-1}. 
\end{equation*}
Choose $L=\lfloor t/l(t)\rfloor$ and $s=t/\ell(t)^{\frac{\al}{2-\al}}$. Then $P\leq \frac{ct^2}{B_t^2l(t)^2}(t/s)^{\frac{2}{\al}-1+\eps}\leq \frac{c}{\ell(t)^{1-\ti\eps}}$.
Thus $P\leq 1/2$ for large enough $t$.
Since $\ell(t)\to\infty$ and $\kappa\to\infty$ as $t\to\infty$, we have that
\begin{equation*}
	\rho_1(2t)\leq \frac{\rho_1(t)}{\ell(t)}+\frac{2l(t)}{t}.
\end{equation*}
Then, for all $t$ such that $\ell(t)>10$,
\begin{equation*}
	\rho_1(2t)-\frac{5\hs l(2t)}{2t}\leq \frac{1}{10}\Big(\rho_1(t)-\frac{5\hs l(t)}{t}\Big),
\end{equation*}
where we used $l(2t)\geq l(t)$ for the increasing function $l$.
The upper bound in (\ref{pribdd2}) follows easily from the last inequality.

{\it Step 3.}
To improve the lower bound to 
\begin{equation}\label{rholow}
	\rho_1(t)\geq \frac{c\hs l(t)}{t},
\end{equation}
we only need to consider the recurrent case, namely when $\al=1$ and $l(\infty)=\infty$, since otherwise the result is trivial.
As in Lemma~\ref{L:nptest} $\Ri_1\leq 0$ and $\Ri_2\leq C$ for some $C>0$, we have that for $s< t$,
\begin{equation}\label{rode5}
	-\frac{\di}{\di t}\rho_1(t) \leq C\rho_1(t-s)^2\sum_xp(x)\P(\hat\tau_x>s)\hs \Big(1+ \frac{1}{B_s\hs \rho_1(t-s)} \Big).
\end{equation}
Recall that $B_t\in\Bi$, there exists $K$ such that $\ell(t)^K\geq l(t)$.
Let $1>\theta>1-1/K$. For each $t$, define $s=s(t)$ from the following equation:
\begin{equation}\label{s}
	s\hs \rho_1(t-s)=l(s)^\theta.
\end{equation}
It follows that $-\rho'_1(t)\leq C\rho_1(t-s)^2/l(s)$ by (\ref{rode5}) and Lemma~\ref{L:jest} and the fact that
\begin{equation}\label{Ri3}
	\frac{1}{B_s\hs\rho_1(t-s)}=\frac{l(s)^{1-\theta}}{\ell(s)}
	\leq \ell(s)^{-1+K(1-\theta)}\underset{t\to0}{\longrightarrow}0.
\end{equation}
By the upper bound of $\rho_1$, for $s=t^\de$ with $\de<1$, the l.h.s.\ of (\ref{s}) is bounded above by $t^{-(1-\de)+\eps}\ll l(s)^{\theta}$.
We thus conclude $ t^{\de}<s<t$.
Moreover, for any $\eps>0$, by the rough upper bound in (\ref{pribdd2}),
\begin{equation}\label{t/s}
	t/s=t\hs \rho_1(t-s)\hs l(s)^{-\theta} \leq C\hs l(t) \hs l(s)^{-\theta} \leq C\hs l(t)^{1-\theta}(t/s)^\eps.
\end{equation}
To simplify the ODE, we have to replace $\rho_1(t-s)^2$ and $l(s)$.
By Lemma~\ref{L:dj} and (\ref{t/s}), for any $\eps>0$,
\begin{equation}\label{d/l}
	1/l(s)-1/l(t)\leq C l(t)^{-1+\eps}\ell(t)^{-1}.
\end{equation}
Furthermore, we use (\ref{drho}) with a change of $t$ to $t-s$, and take $L=2/\rho_1(t-s)$ therein.
It follows that $\frac{Cs}{B_s^2l(s)\rho_1(t-s)}\leq \frac{C}{\ell(s)B_s\rho_1(t-s)}<1$, Thus (\ref{drho}) implies that $\rho_1(t-s)\leq C \rho_1(t)$. 
Therefore, by bounding the r.h.s.\ of (\ref{rode5}) we obtain
\begin{equation*}
	-\frac{\di}{\di t}\rho_1(t) \leq C\rho_1(t)^2/l(t),
\end{equation*}
which yields the lower bound (\ref{rholow}).
\epro

\noi\bpro[of Theorem~\ref{T:den}]
We can now derive the sharp bounds in Theorem~\ref{T:den}.
Recall that for each $t$ we define $s=s(t)$ through \eqref{s}, that is, for $\theta\in(1-1/K,1)$,
\begin{equation}\label{s2}
	s\hs\rho_1(t-s)=l(s)^{\theta}.
\end{equation}
On one hand, the monotonicity of $l$ and the fact that $t>s$ imply 
\begin{equation*}
	s\hs\rho_1(t-s)\leq l(t)^{\theta}.
\end{equation*}
On the other hand, take $\ti s=t/l(\theta)^{\de}$ with $0<\de<1-\theta$.
Using the rough bound $\rho_1(t)\asymp l(t)/t$ from Proposition~\ref{T:priori}, we obtain
\begin{equation*}
	\ti s\hs\rho_1(t-\ti s)\geq t\hs l(t)^{-\de} c\hs l(t)/t =c\hs l(t)^{1-\de},
\end{equation*}
and the r.h.s.\ is larger than $l(t)^{\theta}$ for $t$ sufficient large.
Hence, by monotonicity of $u\mapsto u\rho_1(t-u)$, we conclude that 
\begin{equation*}
	s(t)\leq t/l(t)^{\de}\ll t.
\end{equation*}

Next, recall Lemma~\ref{L:nptest} with $N=2$ and error terms $\Ri_1,\Ri_2$ and $\Ri_3$. 
From the generator calculation $-\big(\rho_1(t)\big)' =\sum_{x\in\Z}p(-x)\E[\xi_t(0)\xi_t(x)]$ as in (\ref{pribdd}), we obtain 
\begin{equation}\label{rhopri}
	-\frac{\di}{\di t}\rho_1(t) = \rho_1(t-s)^2\sum_xp(x)\P(\hat\tau_x>s)\hs \big(1+ \Ri_1+\Ri_2+\Ri_3\big).
\end{equation}
Since the error terms have different expressions for $\al=1$ and $\al<1$, for convenience below we focus on $\al=1$ and $l(\infty)=\infty$, namely the recurrent case. But it is easy to check that all arguments, with slight modifications, work for $\al<1$.
By Lemma~\ref{L:nptest} and the form of $s$,
\begin{equation}\label{Ribdd}
	\Ri_1=O\Big(\frac{1}{l(t)^{1-\theta}}\Big),\quad \Ri_2=O\Big( \frac{1}{\ell(t)^{1-\eps}} \Big),\quand \Ri_3=O\Big( \frac{l(s)^{1-\theta}}{\ell(s)} \Big).
\end{equation}
Applying Lemma~\ref{L:jest} for the time-reversed rate 2 difference process starting from $x$,
\begin{equation}\label{tauxbdd}
	\sum_xp(s)\P(\hat\tau_x>s)=\frac{1}{g_\al(0)\hs l(2t)}\Big(1+O\Big(\frac{1}{\ell(t)^{1-\eps}}\Big)+O\Big(\frac{l(t)^{\eps}}{\ell(t)}\Big)\Big),
\end{equation}
where we approximated $1/l(2s)$ by $1/l(2t)$ using (\ref{d/l}).
Combining (\ref{rode5}), (\ref{Ri3}) and Lemma~\ref{L:jest}, we obtain $-\rho'_1(t)\leq C\rho_1(t-s)^2/l(s)$ for $t$ large. 
Fixing $t$ and the corresponding $s(t)\ll t$, for $u\in [t-s,t]$, the fact that $u\sim t$, together with Proposition~\ref{T:priori}, gives that with $\ti s=s(u)$ as defined through (\ref{s2}),
\begin{equation*}
	-\rho'_1(u)\leq C\rho_1(u-\ti s)^2/l(\ti s)\leq C\rho_1(t-s)^2/l(s).
\end{equation*}
Hence,
\begin{equation}\label{drh}
	\rho_1(t-s)^2-\rho_1(t)^2=-2\int_{t-s}^{t}\rho_1(u)\rho_1'(u)\di u \leq
	C\int_{t-s}^{t}\frac{\rho_1(t-s)\rho_1(t-s)^2}{l(s)}\di u.
\end{equation}
Since now it is known the crude bounds that $\rho_1(t)\asymp l(t)/t$, the inequality $\rho_1(t-s)\leq C\rho_1(t)$ holds. Using this and (\ref{s2}), the r.h.s.\ of (\ref{drh}) is bounded above by $C\rho_1(t)^2/l(s)^{1-\theta}$, and thus
\begin{equation}\label{rhobdd3}
	\rho_1(t-s)^2=\rho_1(t)^2\big( 1+O\big(l(t)^{-(1-\theta)-\eps}\big) \big).
\end{equation}
Combining (\ref{rhopri}), (\ref{Ribdd}) and (\ref{rhobdd3}), choosing $\theta=1-\frac{1}{2K}$, and using $1/\ell(t)\leq 1/l(t)^{1/K}$ and (\ref{t/s}) we achieve that for any $\eps>0$,
\begin{equation}
	-\frac{\di}{\di t}\rho_1(t) = \frac{\rho_1(t)^2}{g_\al(0)\hs l(2t)} \Big(1+ O\Big(\frac{1}{l(t)^{1/2K-\eps}}\Big) \Big),
\end{equation}
which, by Karamata's theorem, gives rise to
\begin{equation}
	\rho_1(t)=\frac{g_\al(0)\hs l(2t)}{t}\Big(1+ O\Big(\frac{1}{l(t)^{1/2K-\eps}}\Big) \Big).
\end{equation}
For the transient case where $\al=1$ and $l(\infty)<\infty$ or $\al<1$, choosing $s$ such that $s\rho_1(t-s)=\ell(s)^{-1/2}$, the same procedure yields that there exists $\eta>0$
\begin{equation}
	\rho_1(t)=\frac{g_\al(0)\hs l(2t)}{t}\Big(1+ O\Big(\frac{1}{\ell(t)^{1/2}\wedge t^{\eta}}\Big) \Big).
\end{equation}
This completes the proof.
\epro

\subsection{Proof of Theorem~\ref{T:npf}}\label{S:thmnpt}
For fixed $\xb\in\Z^N$, the theorem follows from Lemma~\ref{L:nptest} with similar estimates done in Subsection~\ref{S:thmden}. 
More precisely, as before, choose $s$ such that $s\hs \rho_1(t-s)=l(s)^{1-K/2}$ if $\al=1$ and $l(\infty)=\infty$, and choose $s\hs \rho_1(t-s)=\ell(s)^{-1/2}$ otherwise.
The constant $K$ is the one such that $\ell(t)^K\geq l(t)$ as $B_t\in\Bi$.
Then, Theorem~\ref{T:den}, Theorem~\ref{T:noncoa} and Lemma~\ref{L:nptest} imply that there exist constants $c_\xb$ depending on $\xb$ and some $\eta>0$ such that for any $\eps>0$,
\begin{equation}
	\rho_N(\xb,t)=
	\frac{c_\xb\hs g_\al(0)^N}{t^N\hs l(2t)^{\ups_N-N}} 
	\Big( 1+O\Big(\frac{1}{l(t)^{1/2K-\eps} \wedge \ell(t)^{1/2-\eps}\wedge t^{\eta} }\Big) \Big).
\end{equation}

\appendix

\section{Negative correlation}
The negative correlation of coalescing random walks, first proved in \cite[Lemma~1]{Arr81}, is well-known.
We extend the result to a slightly more general form, using a theorem of Harris (see, e.g., \cite[Chapter~II,~Theorem~2.14]{Lig05}).
For convenience, below we identify $\xi_t$ as the set $\{x:\xi_t(x)=1\}$.
For $A\subseteq\Z$, we write $(\xi^A_t)_{t\geq0}$ as the coalescing random walks starting from every point of the set $\xi_0^A=A$.

\begin{theorem}\label{T:necor}
	Let $(\xi^A_t)_{t\geq0}$ be coalescing random walks starting from some $A\subset\Z$. 
	Let $x_1,\ldots,x_n$ be distinct points of $\Z$.
	Then for all $t\geq0$,
	\begin{equation}\label{negcor}
		\P\big(\xi^A_t(x_1,\ldots,x_n)=1\big)
		\leq \P\big(\xi^A_t(x_1,\ldots,x_{n-1})=1\big) \P\big(\xi^A_t(x_n)=1\big).
	\end{equation}
\end{theorem}
\bpro
Let $(\zeta^B_t)_{t\geq0}$, starting from $B$, be the dual voter model constructed from $(\xi^A_t)_{t\geq0}$ using the Harris graphical representation.
For $\cup_i C_i \bigcap D=\emptyset$, define a Markov process $X_t:=(\zeta_t^{C_1},\ldots,\zeta_t^{C_{n-1}},\zeta_t^{D})$ with state space $E:=\{(C,D)\in \Z^{n-1}\times\Z\}$.
Let the partial order of $E$ be $(C,D)\prec (C',D')$ iff $C\subseteq C'$ and $D\supseteq D'$.
Then $X_t$ is a monotone process which jumps only among sites that are comparable under $\prec$. 
Define increasing functions on $(E,\prec)$
\begin{equation}
	f(C,D)= \prod_{i=1}^{n-1}1_{\{C_i\cap A\neq\emptyset\}} \quand
	g(C,D)= 1_{\{D\cap A=\emptyset\}}
\end{equation}
The condition (2.15) of \cite[Chapter~II,~Theorem~2.14]{Lig05} is satisfied, and therefore
\begin{equation}
	\P(\zeta^{C_i}_t\cap A\neq\emptyset,\zeta^{D}_t\cap A=\emptyset)\geq \P(\zeta^{C_i}_t\cap A\neq\emptyset) \P(\zeta^{D}_t\cap A=\emptyset).
\end{equation}
Taking complement of the event $\{\zeta^{D}_t\cap A=\emptyset\}$, letting $C_i=\{x_i\}$ and $D=\{x_n\}$, and using the duality $\{\zeta_t^{x_i}\cap A\neq\emptyset\}=\{\xi_t^{A}\cap \{x_i\}\neq\emptyset\}$, we have (\ref{negcor}).
\epro

\bigskip

\noindent
{\bf Acknowledgements.}
We thank Rongfeng Sun for helpful discussions. We are grateful to the anonymous referees for their careful reading of earlier versions of the manuscript and for their valuable comments.
J.~Yu is supported by National Key R\&D Program of China (No.~2021YFA1002700) and NSFC (No.~12101238 and No.~12271010).

\setlength\bibsep{5pt}
\bibliographystyle{alpha}
\bibliography{bibliography}


\end{document}